\magnification=1200

\def\title#1{{\titlefont\noindent #1\bigskip}}

\def\author#1{{\largefont\noindent #1}\medskip}

\def\beginlinemode{\endmode
 \begingroup\obeylines\def\endmode{\par\endgroup}}
\let\endmode=\par

\newbox\theaddress
\def\address{\smallskip\beginlinemode\parindent 0in\getaddress}
{\obeylines
\gdef\getaddress #1 
 #2
 {#1\gdef\addressee{#2}%
   \global\setbox\theaddress=\vbox\bgroup\raggedright%
    \everypar{\hangindent2em}#2
   \def\endaddress{\egroup\endgroup \copy\theaddress \medskip}}}

\def\thanks#1{\footnote{}{\eightpoint #1}}

\long\def\Abstract#1{{\eightpoint\narrower\vskip\baselineskip\noindent
#1\smallskip}}

\def\skipfirstword#1 {}

\def\ir#1{\csname #1\endcsname}

\newdimen\currentht
\newbox\droppedletter
\newdimen\droppedletterwdth
\newdimen\drophtinpts
\newdimen\dropindent

\def\irrnSection#1#2{\edef\tttempcs{\ir{#2}}
\vskip\baselineskip\penalty-3000
{\largefont\bf\noindent \expandafter\skipfirstword\tttempcs. #1}
\vskip6pt}

\def\irSubsection#1#2{\edef\tttempcs{\ir{#2}}
\vskip\baselineskip\penalty-3000
{\bf\noindent \expandafter\skipfirstword\tttempcs. #1}
\vskip6pt}

\def\irSubsubsection#1#2{\edef\tttempcs{\ir{#2}}
\vskip\baselineskip\penalty-3000
{\noindent \expandafter\skipfirstword\tttempcs. #1}
\vskip6pt}

\def\References{\vbox to.25in{\vfil}\noindent{}{\bf References}
\vskip6pt\par}

\def\References{\vskip6pt\noindent{}{\bf References}
\vskip6pt\par}

\def\baselinebreak{\par \ifdim\lastskip<6pt
         \removelastskip\penalty-200\vskip6pt\fi}

\long\def\prclm#1#2#3{\baselinebreak
\noindent{\bf \csname #2\endcsname}:\enspace{\sl #3\par}\baselinebreak}

\def\rem#1#2{\baselinebreak\noindent{\bf \csname #2\endcsname}:}

\def\qed{{$\diamondsuit$}\vskip6pt}

\def\bibitem#1{\par\indent\llap{\rlap{\bf [#1]}\indent}\indent\hangindent
2\parindent\ignorespaces}

\long\def\eatit#1{}

\def\leftheadlinetext{}
\def\rightheadlinetext{}

\def\leftheadline{{\eightrm\folio\hfil \leftheadlinetext\hfil}}
\def\rightheadline{{\eightrm\hfil\rightheadlinetext\hfil\folio}}

\headline={\ifnum\pageno=1\hfil\else
\ifodd\pageno\rightheadline\else\leftheadline\fi\fi}

\def\tenpoint{\def\rm{\fam0\tenrm}
\textfont0=\tenrm \scriptfont0=\sevenrm \scriptscriptfont0=\fiverm
\textfont1=\teni \scriptfont1=\seveni \scriptscriptfont1=\fivei
\def\mit{\fam1} \def\oldstyle{\fam1\teni}
\textfont2=\tensy \scriptfont2=\sevensy \scriptscriptfont2=\fivesy
\def\cal{\fam2}
\textfont3=\tenex \scriptfont3=\tenex \scriptscriptfont3=\tenex
\def\it{\fam\itfam\tenit} 
\textfont\itfam=\tenit
\def\sl{\fam\slfam\tensl} 
\textfont\slfam=\tensl
\def\bf{\fam\bffam\tenbf} 
\textfont\bffam=\tenbf \scriptfont\bffam=\sevenbf
\scriptscriptfont\bffam=\fivebf
\def\tt{\fam\ttfam\tentt} 
\textfont\ttfam=\tentt
\normalbaselineskip=12pt
\setbox\strutbox=\hbox{\vrule height8.5pt depth3.5pt  width0pt}%
\normalbaselines\rm}

\def\eightpoint{\def\rm{\fam0\eightrm}%
\textfont0=\eightrm \scriptfont0=\sixrm \scriptscriptfont0=\fiverm
\textfont1=\eighti \scriptfont1=\sixi \scriptscriptfont1=\fivei
\def\mit{\fam1} \def\oldstyle{\fam1\eighti}%
\textfont2=\eightsy \scriptfont2=\sixsy \scriptscriptfont2=\fivesy
\def\cal{\fam2}%
\textfont3=\tenex \scriptfont3=\tenex \scriptscriptfont3=\tenex
\def\it{\fam\itfam\eightit} 
\textfont\itfam=\eightit
\def\sl{\fam\slfam\eightsl} 
\textfont\slfam=\eightsl
\def\bf{\fam\bffam\eightbf} 
\textfont\bffam=\eightbf \scriptfont\bffam=\sixbf
\scriptscriptfont\bffam=\fivebf
\def\tt{\fam\ttfam\eighttt} 
\textfont\ttfam=\eighttt
\normalbaselineskip=9pt%
\setbox\strutbox=\hbox{\vrule height7pt depth2pt  width0pt}%
\normalbaselines\rm}

\def\largefont{\def\rm{\fam0\largerm}
\textfont0=\largerm \scriptfont0=\largescriptrm \scriptscriptfont0=\tenrm
\textfont1=\largei \scriptfont1=\largescripti \scriptscriptfont1=\teni
\def\mit{\fam1} \def\oldstyle{\fam1\teni}
\textfont2=\largesy 
\def\cal{\fam2}
\def\it{\fam\itfam\largeit} 
\textfont\itfam=\largeit
\def\sl{\fam\slfam\largesl} 
\textfont\slfam=\largesl
\def\bf{\fam\bffam\largebf} 
\textfont\bffam=\largebf 
\scriptscriptfont\bffam=\fivebf
\def\tt{\fam\ttfam\largett} 
\textfont\ttfam=\largett
\normalbaselineskip=17.28pt
\setbox\strutbox=\hbox{\vrule height12.25pt depth5pt  width0pt}%
\normalbaselines\rm}

\def\titlefont{\def\rm{\fam0\titlerm}
\textfont0=\titlerm \scriptfont0=\largescriptrm \scriptscriptfont0=\tenrm
\textfont1=\titlei \scriptfont1=\largescripti \scriptscriptfont1=\teni
\def\mit{\fam1} \def\oldstyle{\fam1\teni}
\textfont2=\titlesy 
\def\cal{\fam2}
\def\it{\fam\itfam\titleit} 
\textfont\itfam=\titleit
\def\sl{\fam\slfam\titlesl} 
\textfont\slfam=\titlesl
\def\bf{\fam\bffam\titlebf} 
\textfont\bffam=\titlebf 
\scriptscriptfont\bffam=\fivebf
\def\tt{\fam\ttfam\titlett} 
\textfont\ttfam=\titlett
\normalbaselineskip=24.8832pt
\setbox\strutbox=\hbox{\vrule height12.25pt depth5pt  width0pt}%
\normalbaselines\rm}

\nopagenumbers

\font\eightrm=cmr8
\font\eighti=cmmi8
\font\eightsy=cmsy8
\font\eightbf=cmbx8
\font\eighttt=cmtt8
\font\eightit=cmti8
\font\eightsl=cmsl8
\font\sixrm=cmr6
\font\sixi=cmmi6
\font\sixsy=cmsy6
\font\sixbf=cmbx6

\font\largerm=cmr12 at 17.28pt
\font\largei=cmmi12 at 17.28pt
\font\largescriptrm=cmr12 at 14.4pt
\font\largescripti=cmmi12 at 14.4pt
\font\largesy=cmsy10 at 17.28pt
\font\largebf=cmbx12 at 17.28pt
\font\largett=cmtt12 at 17.28pt
\font\largeit=cmti12 at 17.28pt
\font\largesl=cmsl12 at 17.28pt

\font\titlerm=cmr12 at 24.8832pt
\font\titlei=cmmi12 at 24.8832pt
\font\titlesy=cmsy10 at 24.8832pt
\font\titlebf=cmbx12 at 24.8832pt
\font\titlett=cmtt12 at 24.8832pt
\font\titleit=cmti12 at 24.8832pt
\font\titlesl=cmsl12 at 24.8832pt

\tenpoint



\def\manyby{\hbox to.75in{\hrulefill}}
\hsize 5.41667in 
\vsize 7.5in

\def\manyby{\hbox to.75in{\hrulefill}}

\tolerance 3000
\hbadness 3000

\def\item#1{\par\indent\indent\llap{\rlap{#1}\indent}\hangindent
2\parindent\ignorespaces}

\def\itemitem#1{\par\indent\indent
\indent\llap{\rlap{#1}\indent}\hangindent
3\parindent\ignorespaces}

\def\refAH{AH1}
\def\refAHb{AH2}
\def\refBal{B}
\def\refCam{Cam}
\def\refCatb{Cat1}
\def\refCat{Cat2}
\def\refCMa{CM1}
\def\refCMb{CM2}
\def\refDGM{DGM}
\def\refDu{Dub}
\def\refDuVa{DuV1}
\def\refDuVb{DuV2}
\def\refEC{EC}
\def\refE{Ev}
\def\refFL{FL}
\def\refFione{F1}
\def\refFitwo{F2}
\def\refFHH{FHH}
\def\refGGR{GGR}
\def\refGiThesis{Gi1}
\def\refGi{Gi2}
\def\refGib{Gi3}
\def\refHHF{HHF}
\def\refvanc{H1}
\def\refsundance{H2}
\def\refravello{H3}
\def\refpos{H4}
\def\refanti{H5}
\def\reffatpts{H6}
\def\refigp{H7}
\def\refseven{H8}
\def\refnagprob{H9}
\def\refHR{HR}
\def\refHi{Hi}
\def\refHo{Ho}
\def\refId{Id}
\def\refKa{Ka}
\def\refL{L}
\def\refMAC{MAC}
\def\refmig{Mi}
\def\refNone{N1}
\def\refNtwo{N2}
\def\refP{P}
\def\refRoe{R1}
\def\refRoeb{R2}
\def\refS{Seg}
\def\refSe{Sei}
\def\refX{Xu}

\def\C#1{\hbox{$\cal #1$}}

\def\pr#1{\hbox{{\bf P}${}^{#1}$}}
\def\cite#1{[\ir{#1}]}

\def\leftheadlinetext{Brian Harbourne}
\def\rightheadlinetext{Survey of Fat Points on \pr2}

\title{Problems and Progress:}
\vskip-\baselineskip
\title{A survey on fat points in \pr2}

\author{Brian Harbourne}

\address
Department of Mathematics and Statistics
University of Nebraska-Lincoln
Lincoln, NE 68588-0323
email: bharbour@math.unl.edu
WEB: {\tt http://www.math.unl.edu/$\sim$bharbour/}
\smallskip
November 19, 2000\endaddress
\vskip-\baselineskip

\thanks{\vskip -6pt
\noindent This work benefitted from a National Science Foundation grant.
I thank the organizers of the 
International Workshop on Fat Points, February 9-12, 2000, Naples, Italy
for having me participate in such an enjoyable conference.
This paper is based on my talk at the conference, with additional
developments resulting from conversations
with H. Schenck and J. Ro\'e. I also thank A. Geramita and 
M. Johnson for their comments and D. Jaffe
for the use of the laptop on which 
the scripts included here mostly were tested.
\smallskip
\noindent 1980 {\it Mathematics Subject Classification. } 
Primary 13P10, 14C99. 
Secondary 13D02, 13H15.
\smallskip
\noindent {\it Key words and phrases. }  Hilbert function,
resolution, fat points, quasi-uniform, blow up.\smallskip}

\vskip\baselineskip
\Abstract{Abstract: This paper, which expands on a talk
given at the International Workshop on Fat Points, February 9-12, 2000, 
in Naples, Italy, surveys problems and progress on
certain problems involving numerical characters
for ideals $I(Z)$ defining fat points subschemes 
$Z=m_1p_1+\cdots+m_np_n\subset\pr2$ for general points $p_i$.
In addition to presenting some new results,
a collection of MACAULAY 2 scripts
for computing actual or conjectured values of (or bounds on) these
characters is included. 
One such script, {\tt findres}, for example, computes the
syzygy modules in a minimal free resolution of the ideal $I(Z)$
for any such $Z$ with $n\le 8$; since {\tt findres} does not rely on a
Gr\"obner basis calculation, it is much faster than routines that do.}
\vskip\baselineskip

\irrnSection{Introduction}{intro}
This paper surveys work on certain
problems involving fat points subschemes of \pr2.
To encourage experimentation, I have included
a number of MACAULAY 2
scripts for doing explicit calculations.
To simplify using them, I've included them
in the \TeX file for this paper in a verbatim
listing, without any intervening \TeX\ control sequences. 
Thus if you have (or obtain, from, say,
{\tt http://www.math.unl.edu/$\sim$bharbour/Survey.tex}) 
the \TeX listing for this 
paper, you can simply copy the lines for the necessary MACAULAY 2
scripts from this paper directly into MACAULAY, without 
any additional editing.

Although the most general definition of a fat points
subscheme involves the notion of infinitely near 
points (see \cite{reffatpts}), it is simpler here to define a 
{\it fat points subscheme\/} of \pr2 to be a subscheme $Z$ defined 
by a homogeneous ideal $I\subset R$ of the form 
$I(p_1)^{m_1}\cap\cdots\cap I(p_n)^{m_n}$,
where $p_1,\ldots,p_n$ are distinct points of \pr2,
$m_1,\ldots,m_n$ are nonnegative integers and 
$R=k[\pr2]$ is the homogeneous coordinate ring
of \pr2 (i.e., a polynomial ring in 3 variables,
$x$, $y$ and $z$, over an algebraically closed field $k$). 
It is convenient to denote
$Z$ by $Z=m_1p_1+\cdots+m_np_n$ and to denote 
$I$ by $I(Z)$. 

For another perspective,
a homogeneous polynomial $f\in R$ is in $I(Z)$
if and only if $\hbox{mult}_{p_i}(f)\ge m_i$ for all $i$,
where $\hbox{mult}_{p_i}(f)$ denotes the multiplicity
of $f$ at $p_i$ (this being the least 
$t$ such that if $l_1$ and $l_2$ are linear forms
defining lines which meet at $p_i$ and nowhere else then $f$
is in the $t$th power $(l_1,l_2)^t$ of the
ideal $(l_1,l_2)\in R$). 

It is important for what
follows to note that $I(Z)$ is a homogeneous ideal, 
hence $I(Z)$ is the direct sum of its homogeneous components
$I(Z)_t=R_t\cap I(Z)$ (where for each integer $t$, $R_t$ denotes
the $k$-vector space span of all homogeneous polynomials of 
$R$ of degree $t$).

\irSubsection{Numerical Characters}{numchar}
The work I am interested in here concerns
certain numerical characters of ideals 
$I(m_1p_1+\cdots+m_np_n)\subset R$ which take a 
constant value on some nonempty open subset
of points $(p_1,\ldots,p_n)\in(\pr2)^n$.
Thus we will usually consider fat points subschemes
$Z=m_1p_1+\cdots+m_np_n$ for which the points 
$p_i\in\pr2$ are general. (Saying that
something is true for $Z=m_1p_1+\cdots+m_np_n$
for general points $p_i$, is the same as saying that
it holds for some open subset of points
$(p_1,\ldots,p_n)\in(\pr2)^n$.) In order to establish
a result for general points, one
typically establishes it for some particular special
choice of the points and then argues by semicontinuity.
(To justify using semicontinuity, even for specializations
to infinitely near points, see my 1982 thesis,
the relevant parts of which were published in \cite{refsundance};
alternatively,
for specializations keeping the points distinct, see \cite{refP}.)
Thus we will sometimes consider situations for which the
points $p_i$ are in some special position.

Given a fat points subscheme $Z=m_1p_1+\cdots+m_np_n$ and its 
ideal $I=I(Z)$, among the numerical characters which have seen 
attention by various researchers over the years are the following:
\item{$\bullet$} $\alpha(Z)$, the least degree $t$ such that $I(Z)_t\ne 0$;
\item{$\bullet$} $\beta(Z)$, the least degree $t$ such that the zero locus of 
$I(Z)_t$ is zero dimensional;
\item{$\bullet$} $h_Z$, the Hilbert function of $I(Z)$ (i.e.,
the function whose value $h_Z(t)$ for each degree $t$
is the $k$-vector space dimension of $I(Z)_t$);
\item{$\bullet$} $\tau(Z)$, the least degree $t\ge 0$ such that 
$h_Z(t)=P_Z(t)$, where $P_Z$ is the Hilbert polynomial of $Z$
(which is simply $P_Z(s)=(s^2+3s+2-\sum_i m_i(m_i+1))/2$);
\item{$\bullet$} $\nu_t(Z)$, the number of generators of $I(Z)$
in degree $t$ in any minimal set of homogeneous generators.

\noindent (In cases where it is understood which $Z$ is meant, 
I will sometimes write $\alpha$ or $\beta$, etc.,
for the more explicit but more cumbersome $\alpha(Z)$, etc.)

The most fundamental characters are $\alpha$, $h$ and $\nu_t$.
For example, $h_Z$ immediately determines $\alpha(Z)$ and $\tau(Z)$.
Moreover, if one can compute $h_Z$ for any $Z$ then one can also
determine $\beta$ for any particular $Z$. (This is because
$t<\beta(Z)$ if and only if either $t<\alpha(Z)$, or
$t\ge\alpha(Z)$ and there exists some nonzero
$Y=m_1'p_1+\cdots+m_n'p_n$ with $0\le m_i'\le m_i$
for all $i$ such that $h_Z(t)=h_{Z-Y}(t-\alpha(Y))$.
The idea is that for $\alpha(Z)\le t<\beta(Z)$,
every element of $I(Z)_t$ is divisible
by nonconstant homogeneous polynomials $f$
which define divisors in the fixed locus of the linear system
$I(Z)_t$. Any such $f$ spans $I(Y)_{\alpha(Y)}$
for an appropriate $Y$, as above. Since there are only 
finitely many $m_1'p_1+\cdots+m_n'p_n$ with $0\le m_i'\le m_i$, 
one can in principle check whether any such $Y$ exists,
as long as one can always compute $h$.) Similarly, if one can 
compute $\alpha(Z)$ for any $Z$ then one can also determine 
$h_Z$ for any particular $Z$. (Here's how: to compute $h_Z(t)$ 
for some $t$ and some $Z=m_1p_1+\cdots+m_np_n$, let 
$Z_0=Z$ and for each $j>0$ let $Z_j=Z_0+q_1+\cdots+q_j$,
where $q_1,\ldots,q_j$ are general. Since each
additional point $q_i$ imposes one additional condition on 
forms of degree $t$ up to the point where no forms remain,
we see that $h_Z(t)=i$ where $i$ is the least $j$
such that $\alpha(Z_j)>t$.) 

Knowing $\alpha(Z)$ for a particular $Z$ 
sometimes also means we know $\nu_t(Z)$
for all $t$. Indeed, conjectures about 
the values of $\alpha$, and, for certain $Z$, 
of $\nu_t$, are made below (see \ir{HHFconj}). 
There are examples of $Z$ for which, if
$\alpha(Z)$ is what it is conjectured to be, then
so are all $\nu_t(Z)$. (For example, take $Z=m(p_1+\cdots+p_n)$ 
with general points $p_i$, where $n>9$ is an even square and $m$ 
is sufficiently large; see Example 5.2 and Theorem 2.5,
both of \cite{refHHF}.)
Nonetheless, knowing $\alpha$ in general does not seem 
to be enough to determine $\nu_t$ for all $t$,
but if one knows $\nu_t$ for all $t$, then one can always
compute $\alpha$ and hence all of the other characters listed.
(This is because the least $t$ such that $\nu_t(Z)>0$ 
is $t=\alpha(Z)$; moreover, $\nu_{\alpha(Z)}(Z)=h_Z(\alpha(Z))$.)
Thus the characters $\nu_t$ are perhaps even more fundamental 
than the other characters discussed above.

The characters $\nu_t(Z)$ are also interesting due to their 
connection to minimal free graded resolutions of $I(Z)$.
A minimal free graded resolution of $I(Z)$ is an exact sequence
$0\to F_1(Z)\to F_0\to I(Z)\to 0$ in which $F_1(Z)$ and 
$F_0(Z)$ are free graded $R$-modules. It turns out,
up to isomorphism as graded $R$-modules, that $F_0$
is $\oplus_t R[-t]^{\nu_t(Z)}$ and $F_1(Z)$
is $\oplus_t R[-t]^{s_t(Z)}$, where the characters
$s_t(Z)$ are defined via $\nu_t(Z)-s_t(Z)=\Delta^3h_Z(t)$ \cite{refFHH}.
Here $\Delta$ denotes the difference operator (so
for any function $f:{\bf Z}\to {\bf Z}$, we have
$\Delta f(t)=f(t)-f(t-1)$), and $R[i]^j$ denotes the direct sum of $j$
copies of the module $R$ itself but taken with the grading
defined by $R[i]_t=R_{t+i}$.

\irSubsection{Connection to Geometry}{geomconn}
Additional interest in these characters (and essential techniques
in studying them) comes from their connections
to geometry. Given distinct points $p_1,\ldots,p_n\in\pr2$,
let $\pi:X\to\pr2$ be the birational morphism obtained by
blowing the points up. Thus $\pi$ is the unique morphism
where $X$ is a smooth and irreducible rational surface
such that, away from the points $p_i$, $\pi$ is an isomorphism
and such that for each $i$, $\pi^{-1}(p_i)$ is a smooth
rational curve $E_i$. It is known that the divisor class group
$\hbox{Cl}(X)$ is a free abelian group on the classes
$[E_i]$ of the divisors $E_i$ and on the class $[E_0]$,
where $E_0=\pi^{-1}(L)$, $L$ being any line in \pr2 not passing through any
of the points $p_i$. Thus for any divisor $D$ on $X$ we have $[D]
=\sum_{i=0}^na_i[E_i]$ for some integers $a_i$. 

It will be useful
later to recall the intersection form on $\hbox{Cl}(X)$.
This is a symmetric bilinear form denoted
for elements $[C]$ and $[D]$ of $\hbox{Cl}(X)$
by $[C]\cdot[D]$ and determined
by requiring that $[E_i]\cdot[E_j]$ is 0 if $i\ne j$, 1
if $i=j=0$ and $-1$ if $i=j>0$. If $C$ and $D$ are curves
on $X$ such that $C\cap D$ is finite and transverse, then
$[C]\cdot[D]=|C\cap D|$; i.e., $[C]\cdot[D]$
is just the number of points of intersection
of $C$ with $D$. If $C=D$, it is convenient to denote
$[C]\cdot[D]$ by $[C]^2$.

Now, for a divisor
$D$, let $\C O_X(D)$ denote the associated line bundle. 
Given a fat points subscheme $Z=m_1p_1+\cdots+m_np_n$,
it turns out for all $t$ that $h_Z(t)=h^0(X, \C O_X(F_t(Z)))$, where
$F_t(Z)$ is the divisor $tE_0-(m_1E_1+\cdots+m_nE_n)$
and $h^0(X, \C O_X(F_t(Z)))$ denotes the dimension of the 0th cohomology
group $H^0(X, \C O_X(F_t(Z)))$ of the sheaf $\C O_X(F_t(Z))$ 
(i.e., $h^0(X, \C O_X(F_t(Z)))$ is the dimension
of the space of global sections of $\C O_X(F_t(Z))$).
Thus $\alpha(Z)$ is the least $t$ such that $h^0(X, \C O_X(F_t(Z)))>0$
and $\tau(Z)$ is the least $t\ge 0$ such that 
$h^0(X, \C O_X(F_t(Z)))=P_Z(t)$. Moreover, note that
$P_Z(t)=(F_t(Z)^2-K_X\cdot F_t(Z))/2+1$.
By Riemann-Roch we have
$$h^0(X, \C O_X(F_t(Z)))-h^1(X, \C O_X(F_t(Z)))+
h^2(X, \C O_X(F_t(Z)))=P_Z(t),$$
and by duality we know $h^2(X, \C O_X(F_t(Z)))=0$ for $t>-3$,
so $h^0(X, \C O_X(F_t(Z)))-h^1(X, \C O_X(F_t(Z)))=P_Z(t)$ for all $t\ge 0$.
Thus $h_Z(t)=P_Z(t)+h^1(X, \C O_X(F_t(Z)))=
\hbox{dim} R_t - (\sum_i m_i(m_i+1)/2-h^1(X, \C O_X(F_t(Z))))$.
Now, $I(Z)_t$ is precisely
what is left from $R_t$ after imposing for each $i$
the condition of vanishing at $p_i$
to order at least $m_i$; what the previous
equation is saying is that the number of conditions imposed
is $\sum_i m_i(m_i+1)/2-h^1(X, \C O_X(F_t(Z)))$. 
For all $t$ sufficiently large, $h^1(X, \C O_X(F_t(Z)))=0$
so a total of $\sum_i m_i(m_i+1)/2$ conditions are imposed.
For smaller $t$, $h^1(X, \C O_X(F_t(Z)))$ measures the extent 
to which these $\sum_i m_i(m_i+1)/2$ conditions fail to be 
independent, and we can regard $\tau(Z)$ as the least degree 
in which the conditions imposed become independent. 

Likewise, the characters $\nu_t$ can be understood
from two perspectives. There is a natural map
$\mu_t(Z): I(Z)_t\otimes_k R_1\to  I(Z)_{t+1}$
given by multiplication, and $\nu_{t+1}(Z)$ is just 
the dimension of the cokernel of the map $\mu_t(Z)$.
Corresponding to this map $\mu_t(Z)$ we have in a natural way a map
$\mu(F_t(Z)): H^0(X, \C O_X(F_t(Z)))\otimes_k 
H^0(X, \C O_X(E_0))\to  H^0(X, \C O_X(F_{t+1}(Z)))$, and indeed
$\nu_{t+1}(Z)=\hbox{dim}\,\hbox{cok}(\mu(F_t(Z)))$.

\irrnSection{Resolutions}{quresconj}
Given a fat points subscheme $Z=m_1p_1+\cdots+m_np_n$, much current 
work concerns either computing or bounding one or another of 
the numerical characters cited above. Some of the oldest
such work concerned bounding the characters $\nu_t$.

\irSubsection{Dubreil and Campanella Bounds}{dubcamp}
Dubreil \cite{refDu} obtained two bounds on the minimum number
$\sum_i \nu_i(Z)$ of homogeneous generators of $I(Z)$:

\prclm{Theorem}{dubreil}{Let $Z=m_1p_1+\cdots+m_np_n$
be a fat points subscheme of \pr2 with distinct points $p_i$.
Then $\sum_i \nu_i(Z)\le \alpha(Z)+\beta(Z)-\tau(Z)\le \alpha(Z)+1$.}

\noindent{\bf Sketch of proof}: The inequality
$\sum_i \nu_i(Z)\le \alpha(Z)+1$ follows immediately
from the Hilbert-Burch Theorem. Here is a more elementary
proof. Given $R=k[x,y,z]$, we may assume
that $x$, $y$ and $z$ define 
general lines in \pr2 (which, in particular, do not contain any
of the points $p_i$). It is then easy to see that
the image $J=xI(Z)_t + yI(Z)_t$ of the 
map $xI(Z)_t\oplus yI(Z)_t\to I(Z)_{t+1}$
has dimension $2h_Z(t)-h_Z(t-1)$. Since $J$ has a base point
(all elements of $J$ vanish at the common point of vanishing of
$x$ and $y$), we see for all $t\ge\alpha(Z)$ that 
$xI(Z)_t + yI(Z)_t$ cannot contain $zI(Z)_t$.
Hence for $t\ge \alpha(Z)$ the image of $\mu_t(Z)$
has dimension at least $2h_Z(t)-h_Z(t-1)+1$, hence
$\nu_{t+1}(Z)=\hbox{dim}\,\hbox{cok}(\mu_t(Z))
\le h_Z(t+1)-(2h_Z(t)-h_Z(t-1)+1)=\Delta^2h_Z(t+1)-1$,
while of course $\nu_{t}(Z)=h_Z(t)=\Delta^2h_Z(t)$ for $t=\alpha(Z)$.
Summing for $i=\alpha(Z)$ to any $N$ sufficiently
large so that $\nu_j(Z)=0$ and $h_Z(j)=P_Z(j)$ 
for $j\ge N-1$, we obtain $\sum_i\nu_i(Z) \le 1+\sum_i (\Delta^2h_Z(t)-1) 
= 1+\Delta h_Z(N) - (N-(\alpha(Z)-1)) = P_Z(N)-P_Z(N-1) - N + \alpha(Z)
=N+1 - N + \alpha(Z)=\alpha(Z)+1$.

The foregoing proof is based on an argument given 
by Campanella \cite{refCam}. Using a result of \cite{refGGR},
Campanella there also gives a similar but slightly more refined
bound, $\nu_{t+1}(Z)\le \Delta^2h_Z(t+1)-\epsilon_t$,
where $\epsilon_t$ is 0 for $t<\alpha(Z)$, 1 for $\alpha(Z)\le t<\beta(Z)$
and 2 for $\beta(Z)\le t\le\tau(Z)$. Summing these refined
bounds for $i$ from $\alpha(Z)$ to $\tau(Z)+1$ gives
$\sum_i \nu_i(Z)\le \alpha(Z)+\beta(Z)-\tau(Z)$. (This argument requires
that one knows that $\nu_i(Z)=0$ for $i>\tau(Z)+1$, but this is true and 
well known; see \cite{refDGM}. It is also not hard to see this directly,
at least from the point of view of the surface $X$ obtained by blowing up
the points $p_i$. Let $[E_0],\ldots,[E_n]$ be the corresponding 
basis for $\hbox{Cl}(X)$, as discussed in \ir{intro}. The statement we need
to prove is then that $\mu(F_t(Z))$ is surjective if $t>\tau(Z)$.
But $t>\tau(Z)$ means $t-1\ge\tau(Z)$ and hence $h^1(X,\C O_X(F_{t-1}(Z)))=0$,
so taking global sections of the exact sheaf sequence 
$0\to \C O_X(F_{t-1}(Z))\to \C O_X(F_t(Z))\to \C O_{E_0}\otimes\C O_X(F_t(Z))\to0$,
we see $H^0(X, \C O_X(F_t(Z)))$ surjects onto 
$H^0(E_0,\C O_{E_0}\otimes\C O_X(F_t(Z)))$. But since $E_0$ is isomorphic to
\pr1, we know that $\C O_{E_0}\otimes\C O_X(F_t(Z))$ is isomorphic
to $\C O_{E_0}(t)$ and that $H^0(E_0,\C O_{E_0}(1))\otimes H^0(E_0,\C O_{E_0}(t))
\to H^0(E_0,\C O_{E_0}(t+1))$ is surjective, and hence that
$H^0(E_0,\C O_{E_0}(1))\otimes H^0(E_0,\C O_{E_0}\otimes\C O_X(F_t(Z)))
\to H^0(E_0,\C O_{E_0}(t+1))$ is surjective.
Taking global sections (denoted by $\Gamma$)
of $0\to \C O_X\to\C O_X(E_0)\to \C O_{E_0}(1)\to 0$, tensoring
by $V=H^0(X, \C O_X(F_t(Z)))$ and mapping by multiplication, one obtains
the diagram 
$$\matrix{0&\to&\Gamma_X(\C O_X)\otimes V&\to&\Gamma_X(\C O_X(E_0))\otimes 
V&\to&\Gamma_{E_0}(\C O_{E_0}(1))\otimes V&\to&0\cr
           &   &\downarrow&   &\downarrow&   &\downarrow&   & \cr
0&\to&\Gamma_X(\C O_X(F_t(Z)))&\to&\Gamma_X(\C O_X(F_{t+1}(Z)))&\to
&\Gamma_{E_0}(\C O_{E_0}(t+1))&\to&0\cr}$$
in which the leftmost vertical map is obviously an isomorphism and the
rightmost vertical map as we saw is surjective, so
the snake lemma gives an exact sequence
in which the cokernels of the outer vertical
maps are 0, hence the cokernel $\hbox{cok}(\mu(F_t(Z)))$ of the middle
vertical map also vanishes; i.e., $\nu_t(Z)=0$ for $t>\tau(Z)+1$.)      \qed

In addition to Dubreil's bounds in
\ir{dubreil}, and Campanella's upper bounds
$\nu_{t+1}(Z)\le \Delta^2h_Z(t+1)-\epsilon_t$
mentioned in the proof above, Campanella also gave the lower bound
that $\nu_t(Z)\ge \hbox{max}\{\Delta^3h_Z(t),\epsilon_t'\}$,
with $\epsilon_{\beta(Z)}'=1$ and $\epsilon_t'=0$ otherwise 
(these bounds of course follow from
$\nu_t(Z)-s_t(Z)=\Delta^3h_Z(t)$; Campanella was actually 
working in the more general situation of perfect 
codimension 2 subschemes of any projective space).

Here is a result (a slight restatement of Lemma 4.1, \cite{refigp})
that in many cases turns out to be an improvement
on the bounds above of Dubreil and Campanella, proved 
in a way very similar to the proof of \ir{dubreil} given above.
It underlies many of the results of \cite{refigp},
\cite{refseven}, \cite{refHHF} and \cite{refFHH}.
 
\prclm{Lemma}{mybound}{Let $Z=m_1p_1+\cdots+m_np_n$
be a fat points subscheme of \pr2 with distinct points $p_i$
such that $m_1>0$, and let $Z''=Z-p_1$ and $Z'=Z+p_1$.
Then $\hbox{max}(h_Z(t+1)-3h_Z(t)+h_{Z''}(t-1),0)\le 
\nu_{t+1}(Z)\le h_Z(t+1)-3h_Z(t)+h_{Z''}(t-1)+h_{Z'}(t)$.}

\irSubsection{Exact Results}{exact}
Campanella's and Dubreil's bounds hold for any
$Z=m_1p_1+\cdots+m_np_n\subset\pr2$, not just when the 
points $p_i$ are general. Thus it is not surprising
that the bounds are not always exact. For example,
for $Z=3(p_1+\cdots+p_5)$ with $p_i$ general we have
$\nu_8(Z)=2$ but $\hbox{max}\{\Delta^3h_Z(8),\epsilon_8'\}=1$
while $\Delta^2h_Z(8)-\epsilon_7=3$.

Thus we can try to obtain exact results.
The typical pattern for work on fat points has been
first to obtain results when either the multiplicities 
$m_i$ are small or the number $n$ of points is small,
and this is what we see regarding $\nu$.
In particular, for $Z=p_1+\cdots+p_n$ with $p_i$
general, it is easy to see that $\alpha(Z)$ is the least
$t$ such that $t^2+3t+2>2n$ and that $\tau(Z)$
is the least $t$ such that $t^2+3t+2\ge2n$. 
And, as always, $\nu_t(Z)=0$
unless $\alpha(Z)\le t\le\tau(Z)+1$, with,
as always, $\nu_{\alpha(Z)}(Z)=h_Z(\alpha(Z))$, so for
$Z=p_1+\cdots+p_n$ only 
$\nu_{\alpha(Z)+1}(Z)$ remains to be found, and 
Geramita, Gregory and Roberts \cite{refGGR}
proved for such subschemes of general points of \pr2
with multiplicity 1 
that $\mu_{\alpha(Z)}$ has maximal rank
(i.e., $\mu_{\alpha(Z)}$ is either surjective
or injective, and hence $\nu_{\alpha(Z)+1}(Z)=
\hbox{max}\{h_Z(\alpha(Z)+1)-3h_Z(\alpha(Z)), 0\}$).
Using different methods
Id\`a \cite{refId} has extended this to the case that
$Z=2(p_1+\cdots+p_n)$ with $p_i$ general and $n>9$.

Results have also been obtained for large $m$
if $n$ is small. The first such I am aware of is that of
Catalisano \cite{refCat}, who determines $\nu_t(Z)$
for all $t$ and any $Z=m_1p_1+\cdots+m_np_n$ as long as 
the points $p_i$ lie on a smooth plane conic; in particular, this
handles the case of any $Z$ involving
$n\le 5$ general points. 

To discuss Catalisano's result in more detail, let $n$ be 
any positive integer (not necessarily 5 or less) and 
consider $Z=m_1p_1+\cdots+m_np_n$ for any $n$ distinct points $p_i$. 
If $t\ge\alpha(Z)$, let $f_t$ be a common factor 
of greatest degree for the elements of $I(Z)_t$ 
(i.e., $f_t=0$ defines the fixed divisor of the linear system
of curves given by $I(Z)_t$). Now let $Z'_t$
be $\sum_i(\hbox{max}\{m_i,\hbox{mult}_{p_i}(f_t)\})p_i$; thus $f_t$ spans
$I(Z'_t)_d$, where $d=\alpha(Z'_t)$ is the degree of $f_t$. In fact, $d$ 
and $Z'_t$ can be found without finding $f_t$: $d=\alpha(Z'_t)$ where
among all $Z''=m_1''p_1+\cdots+m_n''p_n$ with 
$0\le m_i''\le m_i$ and $h_Z(t)=h_{Z-Z''}(t-\alpha(Z''))$
we choose $Z'_t$ to be that $Z''$ for which $\sum_i(m_i-m_i'')$
is least.

In Catalisano's situation, the points $p_i$ are assumed to lie
on a smooth conic. By \cite{refvanc}, $h_Z$ and $Z_t'$ were already 
known and easy to compute for points on a conic, and, as noted above, 
it is enough to determine $\nu_t(Z)$ for $t>\alpha(Z)$. 
Catalisano's result, although expressed in her paper
\cite{refCat} rather differently, can now be stated:

\prclm{Theorem}{catalisano}{Let $p_1,\ldots,p_n$ 
be distinct points on a smooth plane conic, let
$Z=m_1p_1+\cdots+m_np_n$ be a fat points subscheme of \pr2
and let $d=\alpha(Z'_t)$, where $Z_t'$ is defined as above.
Then for each $t\ge\alpha(Z)$ we have
$\nu_{t+1}(Z)=h_Z(t+1)-h_{Z-Z_t'}(t+1-d)$.}

For a proof in a slightly more general case 
(the conic need not be smooth, for example,
and the points can be infinitely near),
see \cite{reffatpts}. As an aside, note that
\ir{catalisano} shows that $\nu_{t+1}(Z)>0$
only if $\alpha(Z)-1\le t<\beta(Z)$,
since for $t\ge\beta(Z)$ we have $f_t=1$, 
so $d=0$ and $Z_t'=0$.

Since $n>5$ general points do not lie on a conic,
\ir{catalisano} does not apply for $n>5$ general points.
Nonetheless, the inequality $\nu_{t+1}(Z)\ge h_Z(t+1)-h_{Z-Z'_t}(t+1-d)$
always holds (see Lemma 2.10(c), \cite{reffatpts}), although equality
can fail for $n>5$ general points since
$\nu_{t+1}$ can be positive even if $t\ge\beta(Z)$.
A complete solution for the case of any
$Z=m_1p_1+\cdots+m_np_n$ for $n$ general
points $p_i$ was given for $n=6$ by Fitchett \cite{refFitwo},
for $n=7$ by me \cite{refseven} and finally for 
$n=8$ by Fitchett, me and Holay \cite{refFHH}. 

Given $Z=m_1p_1+\cdots+m_np_n$ with $p_i$ general,
the main result of \cite{refFitwo} is that $\mu(F_t(Z))$
has maximal rank as long as $F_t(Z)$ is nef and
$n\le 6$. (A divisor $D$ is {\it nef\/} 
if $[D]\cdot [H]\ge 0$ for every effective divisor $H$.)
Given a divisor on the blow up $X$ of \pr2 at $n$ general points, 
the main results of \cite{refseven} give an algorithmic reduction 
of the problem of determining the rank of $\mu(F)$
to the case that $F$ is ample, and shows if $n\le7$ that
$\mu(F)$ is surjective as long as $F$ is ample,
thereby solving the problem of resolving $I(Z)$ for $n\le 7$.
(A divisor $D$ is {\it ample\/} 
if $[D]\cdot [H]> 0$ for every effective divisor $H$.)

However, if $n=8$, reduction to the ample case 
is not enough, since examples show that $\mu(F_t(Z))$
can fail to have maximal rank even if $F_t(Z)$ is ample;
take $Z=4(p_1+\cdots+p_7)+p_8$ with $t=11$, for instance
(case (c.ii) of \ir{FHHthm}).
The main result of \cite{refFHH}
boils down to giving a formula in nice cases
together with an explicit algorithmic reduction to the nice cases.
We now give a slightly simplified statement of the main result of
\cite{refHHF}.

Recall that an {\it exceptional curve\/}
on a smooth projective surface $S$ is a smooth
curve $C$ isomorphic to \pr1 such that $[C]^2=-1$
in $\hbox{Cl}(S)$. Assuming $n=8$,
let $\Xi_X$ denote the set of classes of exceptional 
curves on $X$ and for each exceptional curve $C$ define
quantities $\lambda_C$ and $\Lambda_C$ as follows:
For $C=E_i$ for any $i$, let $\lambda_C=\Lambda_C=0$.
Otherwise, let $m_C$ be the maximum of $C\cdot E_1, \ldots, C\cdot E_n$,
define $\Lambda_C$ to be the maximum
of $m_C$ and of $(C \cdot L) - m_C$ and define
$\lambda_C$ to be the minimum of $m_C$ and of $(C \cdot L) - m_C$.
We then have:

\prclm{Theorem}{FHHthm}{Let $X$ be obtained by blowing up
8 general points of \pr2 and let 
$[E_0], [E_1],\ldots,[E_8]$ be the associated 
basis of the divisor class group of $X$. Consider the class
$F=t[E_0]-m_1[E_1]-\cdots-m_8[E_8]$, where $m_1\ge \cdots\ge m_8$.  
\item{(a)} If $F \cdot C \ge \Lambda_C$ 
for all $C\in \Xi_X$, then $\mu(F)$ has maximal rank.
\item{(b)} If $F \cdot C < \lambda_C$ for some $C\in \Xi_X$, 
then $\hbox{ker}(\mu(F))$ and $\hbox{ker}(\mu(F-C))$ have the same dimension.
\item{(c)} If neither case (a) nor case (b) obtains, then either
\itemitem{(i)} $F \cdot (E_0 - E_1 - E_2)=0$, in which case 
$\hbox{cok}(\mu_F)$ has dimension 
$h^1(X,{\cal O}_X(F - (E_0 - E_1))) + 
h^1(X,{\cal O}_X(F - (E_0 - E_2)))$, or
\itemitem{(ii)} $[F]$ is $[3E_0 - E_1 - \cdots - E_7] + 
r[8 E_0 - 3E_1 - \cdots - 3E_7 - E_8]$ for some 
$r \ge 1$ (in which case  $\hbox{dim cok}(\mu(F)) = r$ and 
$\hbox{dim ker}(\mu(F)) = r+1$), or
\itemitem{(iii)} $\mu(F)$ has maximal rank.}

This theorem leads directly to an algorithm
for computing resolutions for fat point subschemes $Z$
involving $n\le8$ general points of \pr2.
The MACAULAY 2 script {\tt findres} included at the 
end of this paper implements this algorithm to compute the values
of $\nu_t(Z)$ and $h_Z(t)$ for all $\alpha(Z)\le t\le\tau(Z)+2$. 
Since it does not rely on Gr\"obner basis computations, it is
in comparison quite fast.

\irSubsection{The Quasi-uniform Resolution Conjecture}{resconj}
What to expect for $n>8$ remains mysterious. Whereas
(as discussed below) by taking into account effects due 
to exceptional curves there results
a reasonable conjecture for $h_Z$ for any $Z$ involving
general points, doing the same for resolutions is harder. 
(For a partial result in this direction, see \cite{refFione},
which at least sharpens the bounds Campanella has given
on $\nu_t(Z)$. Also see Theorem 5.3 of \cite{refigp},
which shows in a certain sense that behavior in the $n>8$ case
is simple asymptotically, and that it is the case
of relatively uniform multiplicities that is
not understood.) In fact, the results of \cite{refFHH} for $n=8$
suggest that taking into account the exceptional curves may 
not be enough. Thus it is still unclear how $\nu_t(Z)$
should be expected to behave in general.
 
If one puts a mild condition on the coefficients $m_i$,
however, there is reason to hope that behavior may be
quite simple. In particular, say that $Z=m_1p_1+\cdots+m_np_n$ 
is {\it quasi-uniform\/} if: the points $p_i$ are
general; $n\ge9$ and $m_1=m_9$; and the coefficients
$m_1\ge m_2\ge \cdots\ge m_n\ge 0$ are nonincreasing. 
We then have the following Quasi-uniform Resolution 
Conjecture (\cite{refHHF}):

\prclm{Conjecture}{HHFconj}{If $Z$ is quasi-uniform, then
$h_Z(t)=\hbox{max}\{P_Z(t),0\}$ and
$\nu_t(Z)=\hbox{max}\{h_Z(t)-3h_Z(t-1),0\}$ 
(or, equivalently, $\nu_t(Z)=\hbox{max}\{0, \Delta^3h_Z(t)\}$)
for all $t\ge0$.}

Assuming this conjecture, we can write down an explicit
expression for the resolution of $I(Z)$ for a
quasi-uniform $Z$, as follows (see \cite{refHHF}):
$$0\to R[-\alpha-2]^d\oplus R[-\alpha-1]^c\to 
R[-\alpha-1]^b\oplus R[-\alpha]^a\to I(Z)\to 0,$$
where $\alpha=\alpha(Z)$,
$a=h_Z(\alpha)$, $b=\hbox{max}\{h_Z(\alpha+1)-3h_Z(\alpha),0\}$,
$c=\hbox{max}\{-h_Z(\alpha+1)+3h_Z(\alpha),0\}$, and
$d=a+b-c-1$.

Most of the evidence for this conjecture currently
is for the uniform case (i.e., $Z$ is quasi-uniform with 
all multiplicities $m_i$ equal to some single $m$)
and mostly even then for small multiplicities.
For example, the conjecture is true (and easy) for $n=9$ (see \cite{reffatpts}), 
and whenever $m\le 2$ (see \cite{refGGR} and \cite{refId}). 
However, at the time that this is written,
there are two situations in which the conjecture is known 
for unbounded multiplicities. The first is when $n$
is a power of 4, and $m\ge(\sqrt{n}-2)/4$, in which case
\ir{HHFconj} is true (in characteristic 0) 
by \cite{refHHF} and \cite{refE}.
For the second (again in characteristic 0),
see \cite{refHR}, a typical example of which is $n=d(d+1)$ and 
$m=d+1$ for any even integer $d>2$. 
In this case, the bounds on $\alpha$ and $\tau$ given 
by the modified unloading method due to Ro\'e and me
(using $r=d(2d+1)/2$; see \ir{modunlbnds} and 
the end of \ir{bndst})
show $\alpha(Z)>\tau(Z)$ for $Z=m(p_1+\cdots+p_n)$, 
from which it follows that $I(Z)$ is generated in
degree $\alpha(Z)$ and hence that the minimal free resolution of $I(Z)$ is
\vskip.05in
\noindent\hbox to\hsize{\hfil$0\to R[-d^2-2d-1]^{d(d+2)}\to 
R[-d^2-2d]^{(d+1)^2}\to I(Z)\to 0,$\hfil}
\vskip.05in
\noindent in agreement with \ir{HHFconj}. 
Additional evidence of various kinds for the conjecture is 
given in \cite{refHHF}.

In order to facilitate checking this conjecture and exploring
the problem of understanding resolutions when $Z$ need not
be quasi-uniform, it is helpful to be able to compute
resolutions directly. Although we are interested in
general points, it is easiest instead to use random
choices of points, with the expectation that this will usually
give points that are general enough. It is possible to implement
such a calculation very simply in MACAULAY. Here is an 
example of such a MACAULAY 2 script, provided to me by Hal Schenck,
for computing resolutions of ideals $I(\sum_im_ip_i)$ for random choices
of points $p_i\in\pr2$:

\vfill\eject

\parindent=.7pt
\def\uncatcodespecials{\def\do##1{\catcode`##1=12 }\dospecials}
\def\setupverbatim{\tt\def\par{\leavevmode\endgraf}\obeylines
\uncatcodespecials\obeyspaces}
{\obeyspaces\global\let =\ }%
\def\beginverbatim{\par\begingroup\setupverbatim\doverbatim}
{\catcode`\|=0 \catcode`\\=12 %
|obeylines|gdef|doverbatim^^M#1\endverbatim{#1|endgroup}}

\beginverbatim

R=ZZ/31991[x_0..x_2]

mixer = (l)->(i:=0;
              b:=ideal (matrix {{1}}**R);
              scan(#l, i->(
                  f:=random(R^1,R^{-1});
                  g:=random(R^1,R^{-1});
                  I:=(ideal (f | g))^(l#i);
                  b=intersect(b,I)));
              print betti res coker gens b;
              b)
--Return the ideal of mixed multiplicity random fatpoints. Input is a list
--with the multiplicities; e.g. mixer({1,2,3}) returns the ideal of
--I(p1)^1 \cap I(p2)^2 \cap I(p3)^3, where pj is a (random) point in P^2,
--and prints the betti numbers of the resolution
--HKS 4/28
\endverbatim
\parindent=20pt
\vskip\baselineskip

By a slight modification, below, this script can be made to 
handle the uniform case (i.e., $n$ random points each 
taken with the same multiplicity $m$):

\parindent=.7pt
\def\uncatcodespecials{\def\do##1{\catcode`##1=12 }\dospecials}
\def\setupverbatim{\tt\def\par{\leavevmode\endgraf}\obeylines
\uncatcodespecials\obeyspaces}
{\obeyspaces\global\let =\ }%
\def\beginverbatim{\par\begingroup\setupverbatim\doverbatim}
{\catcode`\|=0 \catcode`\\=12 %
|obeylines|gdef|doverbatim^^M#1\endverbatim{#1|endgroup}}

\beginverbatim

R=ZZ/31991[x_0..x_2]

unif = (n,m)->(<< n << " points of multiplicity " << m << ":" << endl;
               i:=0;
               b:=ideal (matrix {{1}}**R);
               while i < n do (
                   f:=random(R^1,R^{-1});
                   g:=random(R^1,R^{-1});
                   I:=(ideal (f | g))^(m);
                   b=intersect(b,I);
                   i=i+1);
               print betti res coker gens b;
               b)

--Example: unif({3,2}) returns the ideal of
--I(p1)^2 \cap I(p2)^2 \cap I(p3)^2, where pj is a (random) point in P^2,
--and prints the betti numbers of the resolution

\endverbatim
\parindent=20pt

\irrnSection{Hilbert functions}{SHGH}
Nagata, in connection with his work on Hilbert's 14th problem, 
began an investigation of the Hilbert function
$h_Z$ for fat points subschemes $Z=m_1p_1+\cdots+m_np_n$
with $p_i$ general, although his work was written up
from the point of view of divisors on blow ups of \pr2
(see \cite{refNtwo}).

\irSubsection{Nagata's Work}{nagwork}
In brief, Nagata in \cite{refNtwo} determines $h_Z(t)$ (and 
thus $\alpha(Z)$) for all $t$ for any $Z=m_1p_1+\cdots+m_np_n$
with $p_i$ general as long as $n\le9$.
In \cite{refNone}, he poses the following conjecture,
which remains open unless $n$ is a square, in which
case Nagata verified it:

\prclm{Conjecture}{Nagconj}{Let $Z=m_1p_1+\cdots+m_np_n$ for
$n>9$ general points $p_i\in\pr2$. Then
$\alpha(Z)>(m_1+\cdots+m_n)/\sqrt{n}$.}

Also implicit in \cite{refNtwo} is a lower bound
(see ($*$), \ir{decompbnd}) for the values of the Hilbert function 
$h_Z$ of $I(Z)$. An easier lower bound comes from the
fact, as discussed above, that 
$h_Z(t)=P_Z(t)+h^1(X, \C O_X(F_t(Z)))$. Since
$h^1(X, \C O_X(F_t(Z)))\ge 0$, it of course
follows that $h_Z(t)\ge \hbox{max}\{P_Z(t),0\}$.
However, easy examples show that $h_Z(t)>\hbox{max}\{P_Z(t),0\}$
can sometimes occur; in all such examples 
for which $h_Z(t)$ is known, the difference
$h_Z(t)-\hbox{max}\{P_Z(t),0\}$ has a geometric origin,
being always precisely what one gets by taking into account
exceptional curves. Taking the exceptional curves into
account gives the more refined bound ($*$).

To explain this, let $X$ be obtained by blowing up $n$ distinct
points $p_i$ of \pr2. We have, as discussed above,
the basis $[E_0],\ldots,[E_n]$ of the divisor class
group $\hbox{Cl}(X)$ of $X$. Because we are mostly 
interested in the case of $n$ general points, technical
issues force us to use the following definition.
Let us say that an element $v=\sum_ia_i[E_i]$
of $\hbox{Cl}(X)$ is 
an {\it exceptional class\/} if for general points $p_i$
there is an exceptional curve $C\subset X$
with $v=[C]$. (The problem is that there may be
no nonempty open set $U$ of points $(p_1,\ldots,p_n)\in(\pr2)^n$
for which all exceptional classes are simultaneously
classes of exceptional curves, even though each exceptional
class $v$ is the class of an exceptional curve for some nonempty
open $U_v\subset(\pr2)^n$.)

Nagata \cite{refNtwo} determined the set $\C E(n)$ of exceptional classes.
It turns out that $\C E(0)$ is empty,
$\C E(1)=\{[E_1]\}$, and $\C E(2)=\{[E_1],[E_2],[E_0-E_1-E_2]\}$, while
for $n\ge 3$ the set $\C E(n)$ is the orbit $W_n[E_n]$
with respect to the action of the group $W_n$
of linear transformations on $\hbox{Cl}(X)$ generated by
all permutations of $\{[E_1],\ldots,[E_n]\}$ and (if $n\ge 3$)
by the map $\gamma$ for which $\gamma:[E_i]\mapsto [E_i]$ for $i>3$,
$\gamma:[E_i]\mapsto [E_0]-[E_1]-[E_2]-[E_3]+[E_i]$ for $0<i\le 3$
and $\gamma:[E_0]\mapsto 2[E_0]-[E_1]-[E_2]-[E_3]$. (The map $\gamma$
can be regarded as a reflection
corresponding in an appropriate sense to a quadratic Cremona 
transformation centered at $p_1$, $p_2$ and $p_3$. The fact that
$W_n$ is a reflection group was recognized by Du Val \cite{refDuVb},
and extended and exploited by Looijenga \cite{refL}.)


\irSubsection{A Decomposition and Lower Bound}{decompbnd}
For simplicity, assume $n\ge 3$. This will not be a serious
restriction, since cases $n<3$ are easy to handle ad hoc,
and in any case there are natural inclusions
$\C E(n)\subset \C E(n+1)$ for all $n$, so a 
given value of $n$ subsumes smaller values. Now Let $\Psi$ 
be the subsemigroup of $\hbox{Cl}(X)$ generated by 
$\C E(n)$ and by the anticanonical class $-K_X=[3E_0-E_1-\cdots-E_n]$
of $X$. (With respect to the action of $W_n$ on $\hbox{Cl}(X)$,
$\Psi$ is essentially Tits' cone \cite{refKa}; thus 
there exists a fundamental domain for the action of $W_n$ on $\Psi$.)
For any $F\in\Psi$, it turns out that there
is a unique decomposition $F=H_F+N_F$ with 
(dropping the subscripts) $H,N\in\Psi$
such that $H\cdot v\ge 0$ for every exceptional class $v$,
$H\cdot N=0$, and either $N=0$ or $N=a_1v_1+\cdots+a_rv_r$
for some exceptional classes $v_i$ and integers $a_i\ge 0$,
such that $v_i\cdot v_j=0$ for all $i\ne j$.
(It is easy to compute this decomposition. 
By recursively applying $\gamma$ and
permutations in a straightforward 
way, for any $F\in\hbox{Cl}(X)$ one can
find an element $w\in W_n$ such that
either $wF\cdot [E_0]<0$, or $wF\cdot [E_0-E_1]<0$, or 
such that $wF=a_0[E_0]+\sum_{i>0}a_i[E_i]$
with $a_0\ge0$, $a_0+a_1+a_2+a_3\ge0$ and $a_1\le\cdots\le a_n$.
But if either $wF\cdot [E_0]<0$ or $wF\cdot [E_0-E_1]<0$, 
then $F\not\in\Psi$, while otherwise there are two cases. 
Either $0>wF\cdot [E_0-E_1-E_2]=a_0+a_1+a_2$,
in which case $H=w^{-1}[(2a_0+a_1+a_2)E_0-(a_0+a_2)E_1-(a_0+a_1)E_2]$ and 
$N=w^{-1}((-a_1-a_2-a_0)[E_0-E_1-E_2]+\sum_{a_i>0}a_i[E_i])$, 
or $wF\cdot [E_0-E_1-E_2]\ge0$ and
we have $H=w^{-1}(a_0[E_0]+\sum_{a_i<0}a_i[E_i])$
and $N=w^{-1}(\sum_{a_i>0}a_i[E_i])$. An implementation
of this procedure is given by the script {\tt decomp}
provided in this paper.)

It is true (and more or less apparent from \cite{refNtwo})
for general points $p_i$ that if $h^0(X,\C O_X(F))>0$ then
$F\in\Psi$, hence
$F=H+N$ as above, and $h^0(X,\C O_X(F))=h^0(X,\C O_X(H))
\ge (H^2-K\cdot H)/2+1$. For any $F\in\hbox{Cl}(X)$,
define $e(F)$ to be 0 unless $F\in\Psi$, in which case
set $e(F)$ to be the maximum of $1+(H_F^2-K\cdot H_F)/2$ and 0.
We then get the lower bound 
$$\noindent \hbox to\hsize{\hfil $h^0(X,\C O_X(F))\ge e(F)$.\hfil ($*$)}$$
(The script {\tt homcompdim} included at the end of this paper
computes $e(F)$. For $F=d[E_0]-m_1[E_1]-\cdots-m_n[E_n]$,
we have $e(F)=\hbox{\tt homcompdim}(\{d,\{m_1,...,m_n\}\})$,
which, if the multiplicities $m_i$ are nonnegative, 
is also thus a lower bound for the dimension of the homogeneous
component of $I(Z)$ of degree $d$ for $Z=m_1p_1+\cdots+m_np_n$.)

\irSubsection{The SHGH Conjectures}{SHGHconjs}
It follows from Nagata's work (see Theorem 9, \cite{refNtwo}) that in fact
$h^0(X,\C O_X(F))= e(F)$ for $n\le9$ general points. What occurs for $n>9$ is
not known, but I \cite{refvanc} (also see \cite{refravello}), Gimigliano 
\cite{refGiThesis} (also see \cite{refGi}) and
Hirschowitz \cite{refHi} independently gave conjectures for
explicitly computing $h^0(X,\C O_X(F))$ for any $n$.
These conjectures are all equivalent to the following conjecture,
which states that $e(F)$ is the expected value of $h^0(X,\C O_X(F))$:

\prclm{Conjecture}{HGHconj}{Let $X$ be the blow up of $n$ general points 
of \pr2 and let $F\in\hbox{Cl}(X)$. Then $h^0(X,\C O_X(F))=e(F)$.}

It is interesting to compare \ir{HGHconj}
with an earlier conjecture posed by Segre \cite{refS}, giving
a conjectural characterization
of those classes $F$ such that 
$h^0(X,\C O_X(F))>\hbox{max}\{0,(F^2-K\cdot F)/2+1\}$:

\prclm{Conjecture}{Sconj}{Let $X$ be the blow up of $n$ general points 
of \pr2 and let $F\in\hbox{Cl}(X)$. If
$h^0(X,\C O_X(F))>\hbox{max}\{0,(F^2-K\cdot F)/2+1\}$,
then the fixed locus of $|F|$ has a double
component.}

It is easy to show that \ir{HGHconj} 
implies \ir{Sconj}; the fact that \ir{Sconj}
implies \ir{HGHconj} is essentially Theorem 8
of \cite{refNtwo}. Thus I will refer to 
these conjectures (in any of their forms) as
the SHGH Conjecture. Since Nagata's paper is hard 
to read, the equivalence of \ir{HGHconj} 
and \ir{Sconj} was only recently 
recognized and proved by Ciliberto and Miranda. 
Here is a sketch of a proof.

\prclm{Theorem}{equivThm}{\ir{HGHconj} is equivalent to
\ir{Sconj}.}

\noindent{\bf Sketch of proof}: To see that
\ir{HGHconj} implies \ir{Sconj}, assume that 
$h^0(X,\C O_X(F))>\hbox{max}\{0,(F^2-K\cdot F)/2+1\}$
for some $F$. Thus $h^0(X,\C O_X(F))>0$, so $F\in\Psi$
and hence we have a decomposition $F=H+N$, as described above,
with $N=a_1v_1+\cdots+a_rv_r$ for some exceptional classes 
$v_i$ and $a_i\ge 0$. By \ir{HGHconj}, we have 
$0<h^0(X,\C O_X(F))=e(F)=(H^2-K\cdot H)/2+1$, and by substituting
$H+N$ in for $F$ we see that
$(F^2-K\cdot F)/2+1=(H^2-K\cdot H)/2+1$ unless $a_i>1$ for some $i$,
in which case $v_i$ is the class of a curve occurring 
(at least) doubly in the base locus of $|F|$,
proving \ir{Sconj}.

Conversely, assume \ir{Sconj}. Among all $F$ for which
$h^0(X,\C O_X(F))=e(F)$ fails, choose one having
as few fixed components as possible (i.e., for which
the sum of the multiplicities of the fixed components
is minimal). As before we have $F=H+N$, but $N=0$
by minimality (since $h^0(X,\C O_X(H))=h^0(X,\C O_X(F))>e(F)=e(H)$).
Since $F=H$, by construction of $H$ we have $F\cdot E\ge 0$
for every exceptional class $E$. 

Now say some reduced irreducible curve $C$
occurs as a fixed component of $|F|$ with multiplicity at least 2.
Thus $h^0(X,\C O_X(2C))=1$, hence $h^0(X,\C O_X(C))=1$ and
by \ir{Sconj} we have $(C^2-C\cdot K)/2+1=1$ so $C^2=C\cdot K$.
Therefore the genus $g_C$ of $C$ is $(C^2+C\cdot K)/2+1=C^2+1$;
i.e., $C^2\ge -1$. On the other hand $1=h^0(X,\C O_X(2C))
\ge (4C^2-2C\cdot K)/2+1=C^2+1$ so $C^2\le 0$.

If $C^2=-1$, then $g_C=0$, so $C$ is an exceptional curve. From
$0\to \C O_X(F-C)\to \C O_X(F)\to\C O_C(C\cdot F)\to 0$
it follows that $h^1(X,\C O_X(F))=0$ (since
$C\cdot F\ge 0$ implies $h^1(X,\C O_C(C\cdot F))=0$, while
$h^1(X,\C O_X(F-C))=0$ by minimality), which contradicts failure
of $h^0(X,\C O_X(F))=e(F)$. 

If $C^2=0$, then $g_C=1$, so $C$ is an elliptic curve. From
$0\to \C O_X(F-C)\to \C O_X(F)\to\C O_C\otimes \C O_X(F)\to 0$
it follows that $h^1(X,\C O_X(F))=0$ (as long as we see 
$h^1(X,\C O_C\otimes \C O_X(F))=0$, since
$h^1(X,\C O_X(F-C))=0$ by minimality), which again contradicts failure
of $h^0(X,\C O_X(F))=e(F)$. But $C$ (being irreducible of nonnegative
selfintersection) is obviously nef, so $C\cdot F\ge 0$. Since 
$C$ is elliptic, $C\cdot F\ge 0$ guarantees 
$h^1(X,\C O_C\otimes \C O_X(F))=0$ unless the restriction
$\C O_C\otimes \C O_X(F))$ of $\C O_X(F)$ to $C$ is trivial.
But because the points blown up to obtain $X$ are general,
$\C O_C\otimes \C O_X(F)$ cannot be trivial. (In fact, 
up to Cremona transformations, which is to say up 
to the action of $W_n$, we can assume that
$[C]=[3E_0]-[E_1+\cdots+E_9]$. For $\C O_C\otimes \C O_X(F)$
to be trivial we would need $F\cdot C=0$ and, writing $[F]$
as $[dE_0-m_1E_1-\cdots-m_nE_n]$, we would have 
$3d-(m_1+\cdots+m_9)=0$; i.e., $[dE_0-m_1E_1-\cdots-m_9E_9]$
is the class of an effective divisor perpendicular to 
$[3E_0]-[E_1+\cdots+E_9]$, but for general points $p_1,\ldots,p_9$
the only such classes are multiples of $[3E_0]-[E_1+\cdots+E_9]$
itself, for which it is easy to check the restrictions
to $C$ are not, in general, trivial.) \qed

\vfil\eject

\irSubsection{Evidence}{SHGHev}
It is worth mentioning that
it is not hard to show (see \cite{refHHF}) that the \ir{SHGHconj} 
implies that $h_Z(t)=\hbox{max}\{P_Z(t),0\}$
if $Z$ is quasi-uniform, which is part of \ir{HHFconj}
posed above. In particular, if $F=d[E_0]-m[E_1+\cdots+E_n]$
where the $E_i$ are obtained by blowing up $n>9$ general
points of \pr2, then \ir{HGHconj} predicts that
$h^0(X,\C O_X(F))$ equals the maximum of 0 and
$(F^2-K\cdot F)/2+1$. Proving this equality is 
trivial if $m=1$, and was proved 
for $m=2$ by Alexander and Hirschowitz in a series
of papers culminating in \cite{refAH} (it is worth noting that
these papers address \pr{N} for all $N$).
More generally, given any positive integer $M$ and $M\ge m_i>0$ for all $i$,
\cite{refAHb} shows for any $Z=m_1p_1+\cdots+m_np_n$
(in any projective space) that $h_Z(t)=\hbox{max}\{P_Z(t),0\}$
for all $t$ as long as $n$ is sufficiently large compared with $M$.
More explicitly, Ciliberto and Miranda \cite{refCMa}, \cite{refCMb}
have verified the SHGH Conjecture in characteristic 0 
for all $m\le12$ for any $n>9$ (see also \cite{refSe}), 
and Mignon \cite{refmig} has now verified the
\ir{SHGHconj} for all $n>9$ for any
$F=d[E_0]-m_1[E_1]-\cdots-m_n[E_n]$ as long as 
$m_i\le 4$ for all $i$.

Whereas all of the explicit verifications of the SHGH conjecture
described above assume multiplicities at most 12, 
two methods now exist that work for multiplicities
which in some cases can be arbitrarily large; both assume 
that the characteristic is 0. The first is the recent result
of Evain \cite{refE}, which, for example, shows that
$h^0(X,\C O_X(F))$ equals the maximum of 0 and
$(F^2-K\cdot F)/2+1$ for any $F=d[E_0]-m[E_1+\cdots+E_n]$
as long as $X$ is obtained by blowing up $n$ general
points with $n$ being a power of 4. The second  
is the modified unloading method \cite{refHR} jointly due to
me and J. Ro\'e (see \ir{modunlbnds} and the end of \ir{bndst}), 
which gives very tight bounds on $\alpha$ and $\tau$.
With good enough bounds, one can sometimes show 
$\alpha(Z)\ge\tau(Z)$, but anytime one knows $\alpha(Z)\ge\tau(Z)$
it immediately follows that the SHGH Conjecture holds
for $Z$. In fact, there are numerous examples for which
the bounds from \cite{refHR} are good enough to show
$\alpha(Z)\ge\tau(Z)$ and hence that 
the SHGH Conjecture holds, including 
certain infinite families of examples
$Z=m(p_1+\cdots+p_n)$ such as with
$n=d(d+1)$, $m=d+1$ for any even integer $d>2$
(mentioned in \ir{resconj})
or with $n=d^2+2$ and $m=d(d^2+1)+d(d+1)/2$ for any $d>2$
(see Corollary V.2 of \cite{refHR} for these and other examples).

We close this section with the comment
that the script {\tt findhilb} computes the SHGH conjectural values of 
the Hilbert function of $I(m_1p_1+\cdots+m_np_n)$ for general 
points $p_i$. For $n<10$, these values are the actual values.
No separate script for the case of uniform multiplicities
is included since for $Z=m(p_1+\cdots+p_n)$ for $n>9$ general points 
$p_i$, the conjecture is simply that $h_Z(t)$ is the maximum
of 0 and $(t^2+3t+2-nm(m+1))/2$.

\irrnSection{Bounds}{bnds}
Rather than trying to prove the \ir{SHGHconj} directly,
a good deal of work has been directed toward obtaining better
bounds on $\alpha$ and $\tau$. 
The values of $\alpha$ and $\tau$ predicted by the \ir{SHGHconj}
give upper and lower bounds, respectively; in particular,
$\alpha(Z)$ is less than or equal to the least $t$ such that
$e(F_t(Z))>0$, and $\tau(Z)$ is greater than or equal to the 
least $t\ge0$ such that $h_Z(t)=P_Z(t)$.
Thus what is of most interest are lower bounds
on $\alpha$ and upper bounds on $\tau$.
Bounds on $\alpha$ are especially of interest,
since a sufficiently good lower bound on $\alpha$ may
equal the upper bound (and presumed actual value) of
$\alpha$ given by the \ir{SHGHconj}, and, as discussed above,
if $\alpha$ always has its conjectured value
then the full \ir{SHGHconj} is true. 

Unfortunately, such tight bounds are so far fairly rare, 
but there are some, such as $n=d(d+1)$ points taken with multiplicity
$m=d+1$ as discussed above, for which (by 
precisely this method of tight bounds) the Hilbert function and 
resolution are known. For two additional examples, consider
$F=d[E_0]-m[E_1+\cdots+E_n]$ where $X$ is obtained by blowing up
$n$ general points of \pr2 with $n$ being either 16 or 25.
Although Evain's method handles these cases, at least in characteristic
0, an alternate approach is to notice that
the inequality $h^0(X,\C O_X(F))\ge (F^2-K\cdot F)/2+1$ 
with $F=d[E_0]-m[E_1+\cdots+E_n]$ guarantees
that $\alpha\le m\lceil\sqrt{n}\rceil+
\lceil(\lceil\sqrt{n}\rceil-3)/2\rceil$. For
$n=16$ general points this gives $\alpha\le 4m+1$
while for $n=25$ this gives $\alpha\le 5m+1$.
But Nagata's result \cite{refNone} that $\alpha> m\sqrt{n}$ when 
$n$ is a square bigger than 9 now shows that
$\alpha= 4m+1$ for $n=16$ and $\alpha=5m+1$ for $n=25$.
By \cite{refHHF}, $\tau=\alpha$ in these cases, which
determines $h_Z$ for all $t$ (and even the resolution of
$I(Z)$ when $n=16$).

Some of the bounds discussed below are algorithmic in nature,
and hard to give simple explicit formulas or estimates for.
Thus, to compute them, I have included at the end of this paper
two MACAULAY 2 scripts,
{\tt bounds}$(l)$ and {\tt unifbounds}$(l)$; in the former case
$l=\{m_1,\ldots,m_n\}$ (corresponding to taking $n$ general points
with multiplicities $m_1,\ldots,m_n$) while in the latter case
$l=\{n,m\}$ (corresponding to taking $n$ general points
each with multiplicity $m$).


\irSubsection{Bounds on $\alpha$}{bndsa}
By Nagata's work \cite{refNtwo}, the 
exact value of $\alpha(Z)$ is known for any 
$Z=m_1p_1+\cdots+m_np_n$ with $p_i$ general and $n\le9$,
and in such cases can be computed by running the script 
{\tt findalpha} or {\tt uniffindalpha}. 
For $Z=m(p_1+\cdots+p_n)$, $n\le9$, it is easy to be
explicit: $\alpha(Z)=\lceil c_nm\rceil$, where
$c_1=c_2=1$, $c_3=3/2$, $c_4=c_5=2$, $c_6=12/5$,
$c_7=21/8$, $c_8=48/17$ and $c_9=3$.

For $n>9$, {\tt findalpha} or {\tt uniffindalpha} only give upper bounds
for $\alpha$, although the upper bounds given should be, according to 
the SHGH Conjecture, the actual values. Thus most interest is in
finding lower bounds on $\alpha$, and a number of such have been given.
Let $n\ge n'>9$ and consider $Z=m(p_1+\cdots+p_n)$
and $Z'=m(p_1+\cdots+p_{n'})$, where the points $p_i$ are general.
It is easy to see that $\alpha(Z)\ge \alpha(Z')$. Since
Nagata \cite{refNone} proves that $\alpha(Z')>m\sqrt{n'}$ if
$n'$ is a square, it follows (taking $n'=\lfloor\sqrt{n}\rfloor^2$
when $n$ is 16 or more) that $\alpha(Z)>m\lfloor\sqrt{n}\rfloor$.
A complete proof is somewhat tricky; we treat
the slightly weaker inequality
$\alpha(Z)\ge m\lfloor\sqrt{n}\rfloor$ in the next section.

\irSubsubsection{Bounds by testing against nef divisors}{neftest}
The inequality
$\alpha(Z)\ge m\lfloor\sqrt{n}\rfloor$ follows easily (for any $n$)
by specializing $\lfloor\sqrt{n}\rfloor^2$ of the points $p_i$ 
to a smooth plane curve $C'$ of degree $\lfloor\sqrt{n}\rfloor$.
The class $C$ of the proper transform of $C'$ to the blow up $X$
of \pr2 at the points $p_i$ is 
$\lfloor\sqrt{n}\rfloor[E_0]-([E_1+\cdots+E_{\lfloor\sqrt{n}\rfloor^2}])$,
which is nef, but 
$\alpha(Z)[E_0]-m[E_1+\cdots+E_n]$ is (by definition of $\alpha(Z)$)
the class of an effective divisor, so
the intersection $C\cdot (\alpha(Z)[E_0]-m[E_1+\cdots+E_n])=
\alpha(Z)\lfloor\sqrt{n}\rfloor-m\lfloor\sqrt{n}\rfloor^2$
is nonnegative, which gives $\alpha(Z)\ge m\lfloor\sqrt{n}\rfloor$.
More generally, the same argument works for $Z=m_1p_1+\cdots+m_np_n$,
giving $\alpha(Z)\ge (m_1+\cdots+
m_{\lfloor\sqrt{n}\rfloor^2})/\lfloor\sqrt{n}\rfloor$.

Alternatively, by specializing all $n$ points to a curve of degree
$\lceil\sqrt{n}\rceil$, the same argument (using the fact that now
$\lceil\sqrt{n}\rceil[E_0]-[E_1+\cdots+E_n]$
is nef) gives the inequality $\alpha(Z)\ge mn/\lceil\sqrt{n}\rceil$
for $Z=m(p_1+\cdots+p_n)$,
and $\alpha(Z)\ge (m_1+\cdots+m_n)/\lceil\sqrt{n}\rceil$
for $Z=m_1p_1+\cdots+m_np_n$. More generally, we have the 
following extension of the main result of \cite{refnagprob}:

\prclm{Theorem}{bhthm}{Let $Z=m_1p_1+\cdots+m_np_n$
for general points $p_i\in\pr2$ with $n\ge 1$ 
and $m_1\ge \cdots\ge m_n$, and
let $r\le n$ and $d$ be positive integers.
Given nonnegative rational numbers (not all 0)
$a_0\ge a_1\ge \cdots \ge a_n\ge 0$
such that $a_0d^2\ge a_1+\cdots+a_r$
and $ra_0\ge a_1+\cdots+a_n$, then
$\alpha(Z)\ge (\sum_ia_im_i)/(a_0d)$.}

\noindent{\bf Sketch of proof}: Note that by multiplying
by a common denominator, we may assume that each 
$a_i$ is a nonnegative integer.
Consider the class $F=[a_0dE_0-a_1E_1-\cdots-a_nE_n]$
on the surface $X$ obtained by blowing up the points $p_i$.
First, specialize (as in the proof of the main result of 
\cite{refnagprob}) to certain infinitely near points; in particular,
such that $[E_i-E_{i+1}]$ for each $0<i<n$ is 
the class of an effective, irreducible divisor on 
the specialization $X'$ of $X$, and such that 
$d[E_0]-[E_1+\cdots+E_r]$ is the class $C$ of the 
proper transform of a smooth plane curve. Now 
$F$ is nef on $X'$ and hence on $X$. To see this, 
note that: $F\cdot C\ge 0$ since $a_0d^2\ge a_1+\cdots+a_r$;
$F\cdot(E_i-E_{i+1})\ge 0$ for all $i>0$ since $a_i\ge a_{i+1}$;
$F\cdot E_n\ge 0$ since $a_n\ge 0$; and $F$ is a nonnegative
integer sum of the classes $C$, $[E_i-E_{i+1}]$, $i>0$,
and $[E_n]$ since $a_0\ge a_1$ and $ra_0\ge a_1+\cdots+a_n$.
Thus $F$ is a sum of effective classes (in particular,
of $a_0C$ and various multiples of the $[E_i-E_{i+1}]$ and $E_n$), 
each of which it meets nonnegatively; thus $F$ is nef and so $F$ 
meets $\alpha(Z)E_0-m_1E_1-\cdots-m_nE_n$ nonnegatively, 
from which our result follows.   \qed

Finding an optimal bound for a given $Z$
using \ir{bhthm} involves solving a linear
programming problem (note that we may normalize so that $a_0=1$),
not to mention the problem of identifying the best choices of
$r$ and $d$. In case the multiplicities $m_i$ are all equal,
it is not hard to show that optimal solutions (for given $r$ and $d$)
to this linear programming problem are given in
parts (a) and (b) of the following corollary. These need not 
always be optimal if the coefficients are not all equal,
so we consider in parts (c) and (d) some additional possibilities.

\prclm{Corollary}{bhcor}{Let $Z=m_1p_1+\cdots+m_np_n$
for general points $p_i\in\pr2$ with $n\ge 1$ 
and $m_1\ge \cdots\ge m_n\ge0$,
let $r\le n$ and $d$ be positive integers and
let $m$ be the mean of $m_1,\cdots,m_n$.
\item{(a)} If $r^2\ge nd^2$, then $\alpha(Z)\ge mnd/r$.
\item{(b)} If $r^2\le nd^2$, then $\alpha(Z)\ge mr/d$.
\item{(c)} If $d^2\ge r$, then $\alpha(Z)\ge (m_1+\cdots+m_r)/d$.
\item{(d)} Assume $d^2<r$ and let $j$ be an integer, $0\le j\le d^2$.
\itemitem{(i)} If $j=0$, then $\alpha(Z)\ge (m_1+\cdots+m_{d^2})/d$.
\itemitem{(ii)} If $j>0$, let $t=\hbox{min}\{r+(r-d^2)(r-d^2+j)/j,n\}$
and set $m_{t+1}=0$ if $t=n$; then $$\alpha(Z)\ge(1/d)\Big({{(t-\lfloor 
t\rfloor)jm_{t+1}}\over{(r-d^2+1)}}+
\sum_{1\le i\le d^2-j}m_i + \sum_{d^2-j<i\le t}{{jm_i}\over{(r-d^2+j)}}\Big).$$}

\noindent{\bf Sketch of proof}: Each part of the corollary
applies \ir{bhthm} for various values of the $a_i$.
For (a), take $a_0=r$ and $a_1=\cdots=a_n=d^2$.
For (b), take $a_0=n$ and $a_i=r$, $i>0$.
For (c), take $a_i=1$ for $i\le r$ and $a_i=0$ for $i>r$.
For (d)(i), take $a_i=1$ for $i\le d^2$ and $a_i=0$ for $i>d^2$.
For (d)(ii), take $a_i=1$ for $i\le d^2-j$ and $a_i=j/(r-d^2+j)$ 
for $d^2-j<i\le \lfloor t\rfloor$. If $t=n$, then $m_{t+1}=0$
(and so is $t-\lfloor t\rfloor$), but if
$t<n$, then take $a_{t+1}=(t-\lfloor t\rfloor)j/(r-d^2+1)$.

One can formally verify that the values of the $a_i$ given
in (d)(ii) satisfy the necessary conditions to apply \ir{bhthm},
but it may be helpful to briefly discuss how these values come about.
The idea giving rise to the values of $a_i$ in (d)(ii) is to find
extremal sets (one set for each $j$) of values of the $a_i$, 
with the hope that for any given $Z$ one set will be close to an optimal 
solution that might be found by linear programming. By setting $a_0$
equal to 1 (a normalization we clearly can always do), we bound
the values of the other $a_i$ above by 1. Since the multiplicities
$m_i$ are nonincreasing, any optimal solution for the $a_i$ must
also be nonincreasing. Intuitively, we would want to
keep as many of the $a_i$ equal to 1 as possible.
But in order to satisfy $d^2\ge a_1+\cdots+a_r$ we can keep at most
the first $d^2$ of the $a_i$ equal to 1, in which case all of the 
other $a_i$ would have to be 0. Depending on the values of 
the $m_i$, however, we may be better off if we can make enough of the 
other $a_i$ positive. So, given $j$, we leave $a_1,\ldots,a_{d^2-j}$
alone, and spread $a_{d^2-j+1},\cdots,a_{d^2}$, which are 
each 1 to start with, evenly over
$a_{d^2-j+1}$ to $a_r$, which reduces $a_{d^2-j+1},\ldots,a_{d^2}$
from 1 to $j/(r-d^2+j)$, and raises $a_{d^2+1},\ldots,a_r$
from 0 to $j/(r-d^2+j)$, while keeping the condition
$d^2\ge a_1+\cdots+a_r$ satisfied at equality. Now, although this may 
have worsened things (since we may well have reduced
$a_1m_1+\cdots+a_rm_r$), we can hope to more than make up for this
since we can now increase some of the remaining $a_i$ from 0 
(which they were before) to $j/(r-d^2+j)$.
How many of the $a_i$ which we can increase is limited by the
condition $r\ge a_1+\cdots+a_n=d^2+a_{r+1}+\cdots+a_n$; 
moreover, because of fractional effects,
the last $a_i$ which we can manage to increase from 0 might be limited 
to being increased only by a fraction of $j/(r-d^2+j)$, which accounts
for the anomalous behavior of $a_{t+1}$.     \qed

The bounds given in \ir{bhcor} can be computed by running the scripts
{\tt unifbounds} or {\tt bounds}. The script {\tt ezbhalphaD},
which is called by {\tt bounds}, checks all possible $r$, $d$
and $j$ from \ir{bhcor}(d).

\irSubsubsection{Bounds by unloading}{unlbnds}
As an alternative to \ir{bhthm}, we can
use a process that can conveniently be
referred to as {\it unloading}. The 
idea is based on the fact that given a divisor class
$D$ on a surface $X$ and some finite set $S$ of 
classes of effective, irreducible
divisors, if for some $F\in S$ we have $F\cdot D<0$, then
clearly $D$ is the class of an effective divisor
if and only if $D-F$ is. Unloading (in a sense that is slightly
more general than its use in the literature) consists of 
checking $D\cdot F$ for each $F\in S$, and replacing $D$
by $D-F$ whenever $D\cdot F<0$ and continuing with the new $D$.
(For the classical notion of unloading, see pp. 425--438 of vol. 2 of
\cite{refEC}, where it is referred to as {\it scaricamento}, 
or see \cite{refDuVa}.) Eventually, $D$ reduces to a class $D'$ such that 
either $D'$ is obviously not effective
(because, perhaps, $D'\cdot E_0<0$ or $D'\cdot (E_0-E_1)<0$) or
such that $D'\cdot F\ge 0$ for all $F\in S$.

With respect to the specialization used in the proof of 
\ir{bhthm}, we can take $S$ to consist of 
the classes $[E_i-E_{i+1}]$ for $0<i<n$, $[E_n]$ and 
$d[E_0]-[E_1+\cdots+E_r]$, and we look for the largest $t$ such that
$D=[tE_0-(m_1E_1+\cdots+m_nE_n)]$ unloads to 
a class $D'$ with $D'\cdot (E_0-E_1)<0$, in which
case $t+1$ is a lower bound for
$\alpha(Z)$. I have included the script
{\tt bhalpha} to compute the bound obtained via unloading
with respect to any chosen $r$ and $d$.

In the special case that $r\le d^2$, then $d[E_0]-[E_1+\cdots+E_r]$
is nef. It is not hard to then see that the result of the unloading
process is the same as just testing against this nef divisor,
hence, assuming that $m_1\ge \cdots\ge m_n$, we get the bound 
$\alpha(m_1p_1+\cdots+m_np_n)\ge (m_1+\cdots+m_r)/d$.
One can also give a formula for the
result of this unloading process in another extremal case,
$2r\ge n+d^2$. In this case, for $Z=m(p_1+\cdots+p_n)$, we have 
$\alpha(Z)\ge 1+ud+\hbox{min}\{d-1,\lceil\rho/d\rceil-1\}$, where $u\ge0$ and $\rho$
are defined by $mn=ur+\rho$, with $0<\rho\le r$.

The idea of bounding $\alpha$ using unloadings and specializations to 
infinitely near points is due to Ro\'e \cite{refRoe},
who actually uses a sequence of increasingly special
specializations of infinitely near points, applying 
unloading after each specialization. Ro\'e uses a sequence 
of $n-2$ specializations, corresponding to sets $S_i$, $3\le i\le n$, 
of classes of reduced irreducible divisors, where
$S_i=\{[E_n]\}\cup \{[E_j-E_{j+1}]: 1<j<n\}\cup \{[E_1-E_2-\cdots-E_i]\}$.
Starting with $F_t=tE_0-m_1E_1-m_2E_2-\cdots-m_nE_n$,
Ro\'e's algorithm consists of unloading $F_t$ with 
respect to $S_3$ to get $F_t^{(3)}$, then unloading $F_t^{(3)}$ with 
respect to $S_4$ to get $F_t^{(4)}$, etc.,
eventually ending up with 
$F_t^{(n)}=tE_0-m_1^{(n)}E_1-m_2^{(n)}E_2-\cdots-m_n^{(n)}E_n$. 
Ro\'e's bound is then 
$\alpha(Z)\ge m_1^{(n)}$, which comes from the fact that $[F_t^{(n)}]$ 
and hence $[F_t]$ cannot be classes
of effective divisors unless $t\ge m_1^{(n)}$.
This bound can be computed with the scripts
{\tt unifroealpha} and {\tt roealpha}.

Although it is hard to give a simple formula
for the exact value of the result of this method,
an asymptotic analysis by Ro\'e \cite{refRoe} shows that his 
unloading procedure gives a lower bound for $\alpha(Z)$
which is always better than $m(\sqrt{n-1}-\pi/8)$,
for $Z=m(p_1+\cdots+p_n)$ with $n>2$ general points.
It should be noted however that this formula often
substantially understates the result of the full algorithm.

\irSubsubsection{Bounds by a modified unloading}{modunlbnds}
Assume a specialization as in the second paragraph of \ir{unlbnds}; 
we may assume $d[E_0]-[E_1+\cdots+E_r]$ is the class of a smooth curve $C$. 
Unloading $F=[tE_0-m_1E_1-\cdots-m_nE_n]$ with 
respect to the set $S$ of \ir{unlbnds} uses the fact that, if
$C\cdot F<0$, then $F-[C]$ is the class of an effective divisor
if and only if $F$ is. However, the requirement $C\cdot F<0$
can be relaxed, since all we really need is $h^0(C,\C O_C(F))=0$
in order to ensure that $F-[C]$ is the class of an effective divisor
if and only if $F$ is. By joint work with J. Ro\'e \cite{refHR},
using the notion of a flex of a linear series on $C$, one can
show (in characteristic 0) that $h^0(C,\C O_C(F))=0$ if either $t<d$ and
$(t+1)(t+2)/2\le m_1+\cdots+m_r$, or $t>d-3$ and $F\cdot C\le (d-1)(d-2)/2-1$.
(Recall that $C$ is the proper transform of a plane curve $C'$.
The idea is to choose $p_1\in C'$ so that it is not a flex 
for the complete linear series associated to the restrictions to $C$ of 
the divisors occurring during the unloading process.
This is automatic in characteristic 0 as long as $p_1$ is a general point
of $C'$, but in positive characteristics every point of $C'$ may
be a flex for a given, even complete, linear series \cite{refHo}.) 
Using this test in place of the more stringent test $C\cdot F<0$ 
discussed at the beginning of \ir{unlbnds} gives what may be called
the {\it modified unloading\/} procedure. Since this modified procedure 
uses a less stringent test, a larger (or at least as large) degree is needed
to pass the test, so it gives bounds on $\alpha$
which are at least as good as the original unloading procedure.
These new bounds can be computed by running {\tt HRalpha} or {\tt unifHRalpha}.

Although this modified unloading procedure is somewhat difficult to 
analyze in general, in two extremal cases Ro\'e and I 
can derive the following simple bounds for $Z=m(p_1+\cdots+p_n)$, for which
define $u\ge0$, $\rho$ and $s$ by requiring $mn=ur+\rho$, with $0<\rho\le r$,
where $s$ is the largest integer such that $(s+1)(s+2)\le 2\rho$:
\item{$\bullet$} if $2n\ge 2r\ge n+d^2$, 
then $\alpha(Z)\ge s+ud+1$. (The bound given by this formula
can be computed by running {\tt ezunifHRalpha}.)
\item{$\bullet$} if $d(d+1)/2\le r\le \hbox{min}\{n,d^2\}$, 
then $\alpha(n;m)\ge 1+\hbox{min}\{\lfloor(mr+g-1)/d\rfloor,s+ud\}$.
(The bound given by this formula
can be computed by running {\tt ezunifHRalphaB}.)

\irSubsubsection{Bounds using $\Psi$}{psibnds}
The subsemigroup $\Psi\subset \hbox{Cl}(X)$, introduced in 
\ir{decompbnd}, contains the subsemigroup of classes of effective divisors.
Thus, given $Z$, the least $t$ such that $F_t(Z)\in\Psi$ is a lower bound
for $\alpha(Z)$. This sometimes gives an optimal bound. For example,
if $Z=90p_1+80p_2+70p_3+60p_4+50p_5+40(p_6+p_7+p_8)+30p_9+20p_{10}+10p_{11}$,
then the least $t$ such that $F_t(Z)$ is in $\Psi$ is 179,
hence in fact $\alpha(Z)=179$, since $e(F_{179}(Z))>0$.
Finding the least $t$ such that $F_t(Z)\in\Psi$ is somewhat
tedious, so I have provided the script {\tt Psibound} for doing so.

\irSubsubsection{Comparisons}{bndcomps}

For subschemes $Z$ whose multiplicities are not too uniform,
the lower bound on $\alpha(Z)$ given by testing against
$\Psi$ can be the best, as it is for 
$Z=90p_1+80p_2+70p_3+60p_4+50p_5+40(p_6+p_7+p_8)+30p_9+20p_{10}+10p_{11}$
(see \ir{psibnds}).
For example, Ro\'e's method \cite{refRoe} of unloading gives $\alpha(Z)\ge 162$,
and the best result achievable using \ir{bhcor} turns
out to be $\alpha(Z)\ge 173$, whereas testing against $\Psi$
gives $\alpha(Z)\ge 179$ (and hence $\alpha(Z)=179$ as discussed above).

However, if the multiplicities are fairly uniform, 
testing against $\Psi$ does not give a very good bound.
For example, for $Z=m(p_1+\cdots+p_n)$ with $n>9$,
it is easy to see that $F_t(Z)\in\Psi$ for all $t\ge 3m$,
so testing against $\Psi$ gives the bound $\alpha(Z)\ge 3m$.
This compares poorly with bounds via the other methods,
which are typically very close to, but usually less than,
$m\sqrt{n}$. (Currently only the unloading method of \cite{refRoe}
and the modified unloading method, discussed in \ir{unlbnds}
and \ir{modunlbnds}, resp., ever are substantially better than
$m\sqrt{n}$, and even these only when $m$ is not too large 
compared to $n$.)

Thus for uniform subschemes $Z=m(p_1+\cdots+p_n)$ one is better off
using some method other than testing against $\Psi$, such as
testing against nef divisors, as discussed in \ir{neftest}.
In this case one has, for any $d$ and $r$, easy to implement
tests, as given in \ir{bhcor}. By comparison,
the result of unloading with respect to the divisor 
$C=dE_0-(E_1+\cdots+E_r)$, as discussed in \ir{unlbnds},
is, except in certain special cases, 
harder to compute since there is not always a simple formula
for the result. Since one rarely gets something for free,
it is not surprising, for given
$r$ and $d$, that the bounds given by testing against
a nef divisor are never better than those given by unloading.

To see this, let $Z=m(p_1+\cdots+p_n)$, and assume $\alpha(Z)\ge t$
is the bound given by unloading with respect to $C$.
Also, with respect to the same $r$ and $d$, let $F=a_0E_0-(a_1E_1+\cdots+a_nE_n)$
be the nef test class in the proof of \ir{bhthm}. The unloading method
unloads a divisor $D=tE_0-m(E_1+\cdots+E_n)$ to a divisor $D'$
which meets $C$, $E_n$ and $E_i-E_{i+1}$, for all $i$, nonnegatively.
But $D-D'$ is a sum of multiples of these same divisors, which
are all (linearly equivalent to) effective divisors, so each meets $F$ 
nonnegatively. In addition, $F$ is a sum of these same divisors,
each of which $D'$ meets nonnegatively, so $D'\cdot F\ge0$ too. Thus
$D\cdot F\ge 0$, which shows that testing against the nef divisor 
can never rule out the candidate obtained by unloading.

Moreover, if $r^2>nd^2$, unloading can definitely be better.
For example, take $n=22$ and $m=3$. Then the best choice
of $r$ and $d$ with $r^2\ge nd^2$ is $r=19$ and $d=4$,
while the best choice of $r$ and $d$ with $r^2\le nd^2$ 
is $r=14$ and $d=3$. Using \ir{bhcor}(a,b) with 
either choice of $r$ and $d$ gives 
$\alpha\ge 14$, but unloading with respect to
$r=19$ and $d=4$ gives $\alpha\ge 15$. Since \ir{bhcor}
is optimal in this case, we see unloading
sometimes gives a better result than can be obtained by
any use of \ir{bhthm}.

On the other hand, for $Z=m(p_1+\cdots+p_n)$ with  $r^2\le nd^2$ 
and $r\le n$, the bound $\alpha(Z)\ge mr/d$ obtained by 
testing against a nef divisor, cannot be improved by
unloading with respect to $dE_0-(E_1+\cdots+E_r)$, and hence
unloading and testing against a nef divisor give the same
result in these circumstances.
(This is because for unloading to give a better bound,
the class $D=[\lceil mr/d\rceil E_0-m(E_1+\cdots+E_n)]$
would have to unload to something obviously not effective,
but unloading cannot get started unless $D$ meets
$dE_0-E_1-\cdots-E_r$ negatively, which it does not.)
But as the example of the preceding paragraph shows,
if $r^2\le nd^2$, although one cannot do better than $mr/d$
by unloading with respect to $r$ and $d$, one can still hope to
do better than $mr/d$ by unloading using 
some choices $r'$ and $d'$ in place of $r$ and $d$.

Since the modified unloading procedure of \ir{modunlbnds}
uses a less stringent test than does unloading, as in \ir{unlbnds},
with respect to $C=dE_0-(E_1+\cdots+E_r)$ (in the sense that
in order to be allowed to subtract $C$ and continue
the unloading process, for the former the intersection of 
with $C$ can in most cases be as much as $g-1$, 
where $g$ is the genus of $C$, whereas for the latter 
the intersection must be negative), we see that 
bounds obtained via the latter method can never be better than those
obtained by the former. The advantage of the latter method
is that no hypotheses are required on the characteristic.

There is also Ro\'e's unloading method \cite{refRoe}, discussed in \ir{unlbnds}.
As shown in \cite{refnagprob}, for $m$ sufficiently large
compared to $n$, the results of \ir{bhcor} are always better than 
Ro\'e's unloading method. However, when $m$ is not too large compared
with $n$, examples indicate that Ro\'e's method gives the best bounds
currently known. Consider, for instance, two examples 
using modified unloading (\ir{modunlbnds}).
For $n=1000$ and $m=13$, Ro\'e's method gives $\alpha\ge 421$, whereas
modified unloading, using $r=981$ and $d=31$, gives 
$\alpha\ge 424$, and the SHGH conjectural value of $\alpha$ is 426.
For $n=9000$ and $m=13$, things become reversed:
Ro\'e's method gives $\alpha\ge 1274$, while 
modified unloading using
$r=8918$ and $d=94$ gives only $\alpha\ge 1267$;
the SHGH conjectural value of $\alpha$ in this case is 1279.

An interesting feature of these examples is that in both cases
the bounds are better than $\lfloor m\sqrt{n}\rfloor+1$,
conjectured by Nagata (\ir{Nagconj}): For $m=13$ and $n=1000$,
we have $\lfloor m\sqrt{n}\rfloor+1=412$, while for $m=13$ and 
$n=9000$, we have $\lfloor m\sqrt{n}\rfloor+1=1234$.
Indeed, whereas most known lower bounds for $\alpha(Z)$ for $Z=m(p_1+\cdots+p_n)$
are less than $m\sqrt{n}$ (since after all Nagata's conjecture
is still open), the method of \cite{refRoe} and that of 
modified unloading are 
among the few that in certain situations gives bounds that can be substantially
better than $m\sqrt{n}$. In particular, if $m$ is no 
bigger than about $\sqrt{n}$, the method of \cite{refRoe}
consistently (and probably always, although this looks
hard to prove) gives a lower bound that is at least as
big as $m\sqrt{n}$, and gets better as $m$ decreases
until, for $m=1$ it is easy to show that it gives 
the actual value of $\alpha(Z)$. If one chooses $r$ and 
$d$ carefully (depending on $n$), examples indicate
that the modified unloading procedure
does nearly as well as 
the method of \cite{refRoe} when $m$ is small compared to $n$,
and is substantially better for larger $m$. 
The method of \cite{refRoe}, of course, has the advantage
of being characteristic free and does not 
depend on careful choices of other parameters.
The modified unloading method, on the other hand,
sometimes gives a lower bound which is equal to
the SHGH conjectural value (which is known to be an 
upper bound), and thus determines $\alpha$ exactly
(as happens, for example, when $n=d(d+1)$ and $m=d+1$ for
$d>2$ even, as discussed in \ir{resconj}, or $n=38$
with $m=200$, as mentioned in \ir{SHGHev}).

Thus, in terms of getting the best bound
for a given $Z$, the modified unloading method
(at least in characteristic 0) is often the best.
It has, compared with methods (such as \ir{bhcor}) which test 
against nef divisors, the disadvantage of being harder to compute, 
unless special values for $r$ and $d$ are chosen for which a formula 
applies. But since \ir{bhcor} works for essentially any $r$ 
and $d$, sometimes one can do better by applying \ir{bhcor}
than one can by applying the formula of 
\ir{modunlbnds} where one's choices of $r$ and $d$ are more restricted.

This raises the question of which $r$ and $d$ 
give the best result when applying \ir{bhcor}(a, b). In case (a),
$n\ge r$ and $r^2\ge nd^2$ imply $n\ge d^2$ 
(and even $r\ge d^2$), while
in case (b), having $r^2\le nd^2$ and $r\le n$ but trying to
maximize $r/d$ shows that it is enough to consider
values of $d$ with $d\le \lceil\sqrt{n}\rceil$.
In short, in cases (a) and (b), we may as well
only consider $d$ with $d\le \lceil\sqrt{n}\rceil$.
Moreover, given such a $d$, the best choice of $r$
is evidently $\lceil d\sqrt{n}\rceil$ for case (a)
and $\lfloor d\sqrt{n}\rfloor$ for case (b). It is still
(as far as I can see) not easy to tell which $d$ is best 
without checking each $d$ from 1 to $\lceil\sqrt{n}\rceil$,
hence I have included the script {\tt bestrda} for case 
(a), and {\tt bestrdb} for case (b), to do just that. 
Alternatively, $d=\lfloor\sqrt{n}\rfloor$ often seems to be 
a good choice. For this choice of $d$ and the corresponding 
optimal choices of $r$, (a) ends up giving a better bound than 
(b) if $n-d^2$ is even, while (b) is better if $n-d^2$ is odd.

\irSubsection{Bounds on $\tau$}{bndst}
In some ways, $\tau$ is easier to compute than $\alpha$. 
For example, given $Z=m(p_1+\cdots+p_n)$ for $n\ge9$ general 
points, \cite{refHHF} proves by an easy 
specialization argument that 
$$\tau(Z)\le m\lceil\sqrt{n}\rceil
+\lceil(\lceil\sqrt{n}\rceil-3)/2\rceil.$$ 
If $n\ge9$ is a square and
$m>(\sqrt{n}-2)/4$, it follows (see \cite{refHHF}) in fact that
$$\tau(Z)= m\sqrt{n}+\lceil(\sqrt{n}-3)/2\rceil.$$ Thus $\tau$
is known in some situations where $\alpha$ is only conjectured.

Moreover, via an observation of Z. Ran, bounds on $\alpha$ give 
rise to bounds on $\tau$. In particular, given $Z=m(p_1+\cdots+p_n)$
with $p_i$ general, if $\alpha(Z)\ge c_nm$ for all $m$ (where $c_n>0$
depends only on $n$), then 
$$\tau(Z)\le -3+\lceil(m+1)\hbox{max}\{\sqrt{n}, n/c\}\rceil$$
(see Remark 5.2 of \cite{refnagprob}).
Thus, for example, the bounds of \ir{bhcor}(a, b) lead to bounds 
on $\tau$.

It should not be surprising that $\tau$ might be
easier to handle than $\alpha$. Being always able to 
compute $\alpha(Z)$ is equivalent to being always able to compute 
$h_Z$ and hence $\tau(Z)$, while the reverse does not seem
to be true. Moreover, arguments typically involve
specializations. One can hope to compute $\tau$ exactly
using a specialization that drops $\alpha$ (and thereby gives
us something to work with) while leaving $\tau$ unchanged,
but this of course will not work to compute $\alpha$,
only to give a lower bound.

The scripts {\tt findtau} and {\tt uniffindtau} give lower bounds
for $\tau$ which via the SHGH Conjecture are
expected to be the actual values. Thus most 
interest is in finding upper bounds on $\tau$, and
indeed, quite a few upper bounds have been given,
both on \pr2 and in higher dimensions (see,
for example, \cite{refFL}, for various results and
additional references).

Given $Z=m_1p_1+\cdots+m_np_n$, bounding $\tau(Z)$ is mostly of interest
for $n>9$ since for $n\le 9$, for any disposition of the points,
the Hilbert function of $I(Z)$ (and hence $\tau(Z)$)
is known (see \cite{refpos} for $n\le 8$ or \cite{refanti}).
For $n>9$, the results of \cite{refanti} also allow one
to compute $\tau(Z)$ exactly, if the points $p_i$ lie on
a plane cubic. If the points $p_i$ are general, and $t$ is the value 
of $\tau(Z')$ (computed via \cite{refanti}) for some specialization
$Z'$ of the points $p_i$ to a plane cubic, then by semicontinuity
$\tau(Z)\le t$. For $Z=m(p_1+\cdots+p_n)$ with $n>9$, this gives
the bound $$\tau(Z)\le mn/3.$$ This bound is similar in concept to but better
than a bound given by Segre \cite{refS}, obtained by specializing to a 
conic, which for $Z=m(p_1+\cdots+p_n)$ with $n>9$ gives only
$$\tau(Z)\le mn/2.$$

Improved bounds for $\tau(Z)$ for $Z=m_1p_1+\cdots+m_np_n$
with $p_i$ general are given by Catalisano \cite{refCatb}, 
Gimigliano \cite{refGib} and Hirschowitz \cite{refHi}.
For $Z=m_1p_1+\cdots+m_np_n$ with
general points $p_i$ and $m_1\ge\cdots\ge m_n\ge 0$,
Gimigliano's result is that 
$$\tau(Z)\le m_1+\cdots+m_d$$ as 
long as $d(d+3)/2\ge n$, while Hirschowitz's result is that 
$\tau(Z)\le d$ if 
$$\lceil (d+3)/2\rceil\lceil (d+2)/2\rceil>\sum_im_i(m_i+1)/2.$$
Catalisano's result is somewhat complicated, but generalizes
and often improves Gimigliano's. For $Z=m(p_1+\cdots+p_n)$ with 
$n>9$ these all show that
$\tau(Z)$ is at most approximately $m\sqrt{2n}$.
For $n$ sufficiently large, this clearly is better than 
$\tau(Z)\le mn/3$. 

The bound $\tau(Z)\le m\lceil\sqrt{n}\rceil
+\lceil(\lceil\sqrt{n}\rceil-3)/2\rceil$, mentioned 
above (\cite{refHHF}), results from specializing
$n>9$ points to a smooth curve of degree $\lceil\sqrt{n}\rceil$.
Two other bounds which are also on the order of $m\sqrt{n}$
are Ballico's \cite{refBal} for which 
$\tau(Z)\le d$ if $$d(d+3)-nm(m+1)\ge 2d(m-1)-2$$
(but note that $\tau(Z)\le m\lceil\sqrt{n}\rceil
+\lceil(\lceil\sqrt{n}\rceil-3)/2\rceil$ is better for any given
$n$ if $m$ is large enough)
and Xu's for which $\tau(Z)\le d$ if $$3(d+3)>(m+1)\sqrt{10n}$$
(although $\tau(Z)\le m\lceil\sqrt{n}\rceil
+\lceil(\lceil\sqrt{n}\rceil-3)/2\rceil$ is better if $n$
is sufficiently large).

By employing a sequence of specializations
to infinitely near points similar to what he did for bounding
$\alpha$, Ro\'e \cite{refRoeb} obtains an upper bound 
on $\tau$. The method
applies for any $Z=m_1p_1+\cdots+m_np_n$, with $p_i$ general
and $n\ge 2$. For $Z=m(p_1+\cdots+p_n)$, \cite{refRoeb} denotes
this upper bound by $d_1(m,n)$ and proves
$$d_1(m,n)+1\le m(n/(n-1))(\Pi_{i=2}^{n-1}((n-1+i^2)/(n-1+i^2-i)))
\le (m+1)(\sqrt{n+1.9}+\pi/8).$$
The bound $\tau(Z)\le (m+1)(\sqrt{n+1.9}+\pi/8)-1$ compares very well
with the bound $\tau(Z)\le m\lceil\sqrt{n}\rceil
+\lceil(\sqrt{n}-3)/2\rceil$: the former is better
for approximately 60\% of the values of $n$ between
any two successive squares.

Given a curve $C$, the idea of Ro\'e's algorithm is that
for any $F$, by taking cohomology of
$0\to \C O_X(F-C)\to \C O_X(F)\to \C O_C\otimes\C O_X(F)\to 0$,
we have $h^1(X,\C O_X(F))=0$ if $h^1(X,\C O_X(F-C))=0$ and
$h^1(X,\C O_C\otimes\C O_X(F))=0$. In Ro\'e's case, $C$
is always rational so $h^1(X,\C O_C\otimes\C O_X(F))=0$
is guaranteed if $F\cdot C>-2$, and he handles
$h^1(X,\C O_X(F-C))=0$ by induction. 

In somewhat more detail, start with $Z=m_1p_1+\cdots+m_np_n$
with $p_i$ general, and so we may assume $m_i\ge m_{i+1}\ge 0$ 
for all $i$. We have the corresponding divisor class
$F=[tE_0-m_1E_1-\cdots-m_nE_n]$ where $t$ is as yet undetermined.
Now specialize so that each element of 
$S=\{[E_i-E_{i+1}]: 1<i<n\}\cup\{[E_n]\}$
and $[E_1-E_2]$ is the class of a reduced, irreducible divisor.
Now, $F\cdot [E_1-E_2]\ge-1$ is certainly true
to start with (in fact, we have $F\cdot [E_1-E_2]\ge0$).
If $F\cdot [E_1-E_2-E_3]\ge-1$, then fine, but otherwise
replace $F$ by $F-[E_1-E_2]$ and unload the result
with respect to $S$, and continue replacing and unloading
in the same way until $F\cdot [E_1-E_2-E_3]\ge-1$. Note that
throughout this sequence of operations we have $F\cdot [E_1-E_2]\ge-1$, 
so (taking $[C]=[E_1-E_2]$) we have $h^1(X,\C O_C\otimes\C O_X(F))=0$.
Also, unloading involves a succession of replacements of 
$F$ by $F-[E]$, where $[E]$ is always either $[E_i-E_{i+1}]$ for some
$i$ or $[E_n]$, and can always be carried out 
in such a way that at each step we have $F\cdot [E]>-2$. 
Thus we always have $h^1(E,\C O_E\otimes\C O_X(F))=0$, 
where $E$ is a curve whose class is,
at various times, $[E_i-E_{i+1}]$ for some $i$ or $[E_n]$.

So eventually $F$ turns into a class
for which $F\cdot [E_1-E_2-E_3]\ge -1$,
$F\cdot [E_i-E_{i+1}]\ge 0$ for all $i$ and $F\cdot E_n\ge 0$. 
We now further specialize
so that $[E_1-E_2-E_3]$ is the class of an irreducible divisor,
and keep replacing $F$ by $F-[E_1-E_2-E_3]$, unloading with 
respect to $S$ after each replacement, as long as 
$F\cdot [E_1-E_2-E_3-E_4]<-1$. We continue in this way, 
specializing successively so that each 
$[E_1-E_2-\cdots-E_i]$ in turn becomes 
the class of an irreducible divisor,
and replacing $F$ by $F-[E_1-\cdots-E_i]$
and unloading with respect to $S$ after each replacement, as long as 
$F\cdot [E_1-E_2-\cdots-E_{i+1}]<-1$. Eventually
we end up with a class $F'=[tE_0-m_1'E_1-\cdots-m_n'E_n]$ with
$F'\cdot [E_i-E_{i+1}]\ge0$ for all $i$, $F\cdot E_n'\ge0$ and
$F'\cdot [E_1-E_2-\cdots-E_n]\ge-1$.
By construction, $h^1(X,\C O_X(F))=0$ for our original class $F$
if $h^1(X,\C O_X(F'))=0$, but it turns out in
the specialization we end up with that $h^1(X,\C O_X(F'))=0$
if $t\ge m_1'+m_2'-1$. Thus Ro\'e's bound is 
$\tau(Z)\le m_1'+m_2'-1$.

By combining (in characteristic 0) 
the approaches of \cite{refRoeb}, \cite{refnagprob} and \cite{refHHF},
similar to what is done in \ir{modunlbnds},
Ro\'e and I \cite{refHR} obtain another bound on $\tau$. 
The method uses a single specialization
in which the same set $S$ as above consists of classes
of irreducible divisors, but instead of $[E_1-E_2-\cdots-E_{i+1}]$
being the class of an irreducible divisor $D$, 
$[dE_0-(E_1+\cdots+E_r)]$ is, for some $d$ and $r\le n$.
The idea is to start with some class 
$F=[tE_0-m_1E_1-\cdots-m_nE_n]$ with $m_i\ge m_{i+1}\ge 0$ 
for all $i$. We want to choose $t$ to be large enough to start with
so that we can keep subtracting $[dE_0-(E_1+\cdots+E_r)]$
and unloading with respect to $S$ until we eventually 
obtain a class $F'=t'[E_0]$ for some $t'$, while along the 
way always keeping $h^1(D,\C O_D\otimes \C O_X(F))=0$.
The latter is guaranteed (in characteristic 0) if
both $F\cdot E_0\ge d-2$ and 
$F\cdot D\ge g-1$, where $g=(d-1)(d-2)/2$ is the genus of $D$.

The output of the algorithm of the previous paragraph
is easy but tedious to compute in any given case; to get a 
nice formula we seem to need to choose $r$ and $d$ carefully.
For example, let $Z=m(p_1+\cdots+p_n)$ with $p_i$
general in characteristic 0. Assume $r\le n$ and define
$u\ge0$ and $0<\rho\le r$ via $mn=ur+\rho$.
If $2r\ge n+d^2$ (such as is the case for
$d=\lfloor\sqrt{n}\rfloor$ and 
$r=\lceil d\sqrt{n}\rceil$), the algorithm gives 
$$\tau(Z)\le \hbox{max}\{\lceil(mr+g-1)/d\rceil,(u+1)d-2\},$$
while if $r\le d^2$, then the algorithm gives
$$\tau(Z)\le \hbox{max}\{\lceil(\rho+g-1)/d\rceil+ud,(u+1)d-2\}.$$
Using $d=\lceil\sqrt{n}\rceil$ and $r=n$, the latter formula
gives a bound which is always at least as good as that mentioned above
from \cite{refHHF}. And when $m$ is sufficiently large,
the former formula becomes
$\tau(Z)\le \lceil mr/d+(d-3)/2\rceil$, which for a given $n$
with $m$ sufficiently large, gives a better 
bound than the bound $d_1(m,n)$ given in \cite{refRoeb}.
(To justify this claim, note that
by a method similar to how \cite{refRoeb} shows that 
$d_1(m,n)\le -1+m(n/(n-1))(\Pi_{i=2}^{n-1}((n-1+i^2)/(n-1+i^2-i)))$
one can show that
$m(n/(n-1))(\Pi_{i=2}^{n-1}((n-1+i^2)/(n-1+i^2-i)))-
\sum_{i=3}^nn/(i(n-1)+i(i-1)(i-2))\le d_1(m,n)$.
But $(n/(n-1))(\Pi_{i=2}^{n-1}((n-1+i^2)/(n-1+i^2-i)))
\ge n/(\sqrt{n-1}-\pi/8+1/\sqrt{n-1})$; 
see the proof of Proposition 4.2 of \cite{refnagprob}.
The claim now follows for $m$ large enough from the fact that
$n/(\sqrt{n-1}-\pi/8+1/\sqrt{n-1})>r/d$ for
$d=\lfloor\sqrt{n}\rfloor$ and 
$r=\lceil d\sqrt{n}\rceil$ when $n\ge 10$.)

The formulas 
$$\tau(Z)\le \hbox{max}\{\lceil(mr+g-1)/d\rceil,(u+1)d-2\}$$
and
$$\tau(Z)\le \hbox{max}\{\lceil(\rho+g-1)/d\rceil+ud,(u+1)d-2\}$$
can be evaluated by running
{\tt ezunifHRtau} and {\tt ezunifHRtauB}, respectively. 
Since the algorithm works for any $r\le n$ and $d$, it can sometimes
do better than the formulas, which only work for certain 
values of $r$ and $d$. Thus I have provided scripts {\tt unifHRtau}
and {\tt HRtau} to compute the output of the full algorithm
with respect to any specified choice of $r\le n$ and $d$. 

\irrnSection{Scripts}{scrpts}
We close this survey with a collection of MACAULAY 2 scripts
for computing some of the quantities and bounds discussed above.
This is a verbatim listing: There are no \TeX\ control sequences
interspersed in the text of the scripts in the \TeX file for this paper,
so one can simply copy the text of the scripts
from the \TeX file directly into a file called (say)
{\tt BHscripts}. To run a script, such as {\tt findres} (which computes
a resolution of $I(Z)$ for $Z=m_1p_1+\cdots+m_8p_8$, where the $p_i$ 
are assumed to be general
and each $m_i$ is an integer), start MACAULAY 2 and enter
the command {\tt load "BHscripts"}. Then enter the command
{\tt findres($\{m_1,m_2,m_3,m_4,m_5,m_6,m_6,m_7,m_8\}$)}. 

The required format 
for each script's input parameters are described below,
just before the listing for each script. Individual scripts can be run
without loading the entire file, but many scripts defined below call 
one or more of the others, so be sure to load all scripts
called by the one you wish to run.
\vskip\baselineskip

\eightpoint

\parindent=.7pt
\def\uncatcodespecials{\def\do##1{\catcode`##1=12 }\dospecials}
\def\setupverbatim{\tt\def\par{\leavevmode\endgraf}\obeylines
\uncatcodespecials\obeyspaces}
{\obeyspaces\global\let =\ }%
\def\beginverbatim{\par\begingroup\setupverbatim\doverbatim}
{\catcode`\|=0 \catcode`\\=12 %
|obeylines|gdef|doverbatim^^M#1\endverbatim{#1|endgroup}}

\beginverbatim
-- These routines have been debugged on MACAULAY 2, version 0.8.52
-- Brian Harbourne, October 12, 2000

-- findres: This computes the syzygy modules in any resolution
-- of the saturated homogeneous ideal defining any eight or fewer general 
-- fat points of P2. The hilbert function of the ideal is also found.
-- Call it as findres({m_1,...,m_n}) for n <= eight integers m_i.
-- Note that findres does not rely on Grobner bases, so it is fast by comparison.

findres = (l) -> (
if #l>8 then (
   << "This script works only for up to 8 points." << endl;
   << "Please try again with an input list of at most 8 integers." << endl)
else (
   i:=0;
   myflag2:=0;
   w2:={};
   dd1:=0;
   myker:=0;
   www:={};
   ww:=l;
-- the list l of multiplicities is, for simplicity, extended if need be
-- so that it has 8 elements.
   while(#ww < 8) do ww=join(ww,{0});
   n:=#ww;
   n=n-1;
   ww=zr(ww);  -- zero out negative elements of the list
   a1:=findalpha(ww);  -- find alpha, the least degree t such that I_t \ne 0
   d1:=a1-2;
   tau:=findtau(ww);
   v4:={}; -- list of number of syzygies in each degree t listed in v0
   v3:={}; -- list of dim of coker of \mu_t in each degree t listed in v0
   v2:={}; -- list of dim of ker of \mu_t in each degree t listed in v0
   v1:={}; -- list of Hilbert function values for each degree listed in v0
   v0:={}; -- list of degrees from alpha-2 to tau+2, where tau is the least 
         -- degree such that the fat points impose independent conditions
   while (d1 <= tau+2) do ( -- loop from alpha-2 to tau+2, computing v0, v1 and v2
-- append the current degree d1 to the list of degrees
      v0=join(v0,{d1});
-- append the value of the hilbert function in degree d1
-- to the list of values of the hilbert function
      v1=join(v1,{homcompdim(fundom({d1,ww}))});
-- now compute and append to the list v2 the dimension myker of the 
-- kernel of the map \mu_t : I_t\otimes k[P2]_1 \to I_{t+1} where t=d1 
-- and I is the ideal of the fat points subscheme
      if d1<a1 then (v2=join(v2,{0}));  -- d1<a1 means I_{d1}=0 so myker=0
      if d1>=a1 then (     -- for d1>=a1, compute myker 
         myflag2=0;  
         w2=ww;
         dd1=d1;  
         while(myflag2==0) do ( -- this loop implements the main theorem of [FHH]
                                -- which gives an algorithm for computing myker
            w2=zr(w2);  
            w2=prmt(w2);  
            if homcompdim({dd1,w2})==0 then myflag2=1 else (  
              if dd1*6 - (dot(w2,{3,2,2,2,2,2,2,2})) <= 2 then (  
                 dd1=dd1-6;  
                 w2 = {(w2#0)-3,(w2#1)-2,(w2#2)-2,(w2#3)-2,(w2#4)-2,
                 (w2#5)-2,(w2#6)-2,(w2#7)-2}) else (  
                 if dd1*5 - (dot(w2,{2,2,2,2,2,2,1,1})) <= 1 then (  
                    dd1=dd1-5;  
                    w2 = {(w2#0)-2,(w2#1)-2,(w2#2)-2,(w2#3)-2,(w2#4)-2,
                    (w2#5)-2,(w2#6)-1,(w2#7)-1}) else (  
                    if dd1*4 - (dot(w2,{2,2,2,1,1,1,1,1})) <= 1 then (  
                       dd1=dd1-4;  
                       w2 = {(w2#0)-2,(w2#1)-2,(w2#2)-2,(w2#3)-1,(w2#4)-1,
                       (w2#5)-1,(w2#6)-1,(w2#7)-1}) else (  
                       if dd1*3 - (dot(w2,{2,1,1,1,1,1,1,0})) <= 0 then (  
                          dd1=dd1-3;  
                          w2 = {(w2#0)-2,(w2#1)-1,(w2#2)-1,(w2#3)-1,(w2#4)-1,
                          (w2#5)-1,(w2#6)-1,(w2#7)}) else (  
                          if dd1*2 - (dot(w2,{1,1,1,1,1,0,0,0})) <= 0 then (  
                             dd1=dd1-2;  
                             w2 = {(w2#0)-1,(w2#1)-1,(w2#2)-1,(w2#3)-1,(w2#4)-1,
                             (w2#5),(w2#6),(w2#7)}) else (  
                             if dd1 - (dot(w2,{1,1,0,0,0,0,0,0}))< 0 then (  
                                dd1=dd1-1;  
                                w2 = {(w2#0)-1,(w2#1)-1,(w2#2),(w2#3),(w2#4),
                                (w2#5),(w2#6),(w2#7)}) else (  
                                myflag2=2))))))));
         if myflag2==1 then myker=0 else (  
            if dd1 - (dot(w2,{1,1,0,0,0,0,0,0})) == 0 then   
            myker=homcompdim({dd1-1,{(w2#0)-1,(w2#1),(w2#2),(w2#3),
               (w2#4),(w2#5),(w2#6),(w2#7)}})+  
            homcompdim({dd1-1,{(w2#0),(w2#1)-1,(w2#2),(w2#3),(w2#4),(w2#5),
               (w2#6),(w2#7)}}) else (  
               www={3*(w2#7)+1,3*(w2#7)+1,3*(w2#7)+1,3*(w2#7)+1,
               3*(w2#7)+1,3*(w2#7)+1,3*(w2#7)+1,(w2#7)};  
               if {8*(w2#7)+3,www}=={dd1,w2} then myker=(w2#7)+1 else (
                  if homcompdim({dd1+1,w2})>3*(homcompdim({dd1,w2})) then 
                  myker=0 else myker=3*(homcompdim({dd1,w2}))-homcompdim({dd1+1,w2})))); 
         v2=join(v2,{myker}));  
      d1=d1+1);  
   scan(#v0, i->( -- this scan computes v3 from v2 and v1
         if i<2 then v3=join(v3,{0}) else (  
            if v2#(i-1)>-1 then v3=join(v3,{v1#i-3*(v1#(i-1))+v2#(i-1)}) else (  
               if v1#i-3*(v1#(i-1))>=0 then 
               v3=join(v3,{v1#i-3*(v1#(i-1))}) else v3=join(v3,{0})))));  
   scan(#v0, i->( -- this scan computes v4 from v3 and v1
         if i<3 then v4=join(v4,{0}) else 
         v4=join(v4,{v3#i-v1#i+3*(v1#(i-1))-3*(v1#(i-2))+v1#(i-3)})));    
   << "The output matrix has four columns. Column 1 indicates" << endl;
   << "each degree from alpha-2 (where alpha is the least" << endl;
   << "degree t such that I_t > 0 for the fat points ideal I)" << endl;
   << "to tau+2 (where tau is the least degree t such that the points" << endl;
   << "impose independent conditions in all degrees t or bigger)." << endl;
   << "Column 2 gives the value dim I_t of the Hilbert function" << endl;
   << "in each degree t listed in column 1. The resolution of I" << endl;
   << "is of the form 0 -> F_1 -> F_0 -> I -> 0, where F_1 and" << endl;
   << "F_0 are free S=k[P2] modules. Thus F_0=oplus_t S[-t]^{n_t}" << endl;
   << "and F_1=oplus_t S[-t]^{s_t} for integers s_t and n_t." << endl;
   << "Columns 3 and 4 give the values of n_t and s_t in each degree t" << endl;
   << "listed in column 1 (n_t and s_t are 0 in all other degrees)." << endl;
   transpose(matrix({v0,v1,v3,v4}))))

-- findhilb: computes e(F_t(Z)) for alpha-1<= t <=tau+1, which gives 
-- a lower bound for the SHGH conjectural hilbert function
-- for a fat points subscheme involving general points of P2.
-- Call it as findhilb({m_1,...,m_n}) for integers m_i
-- specifying the multiplicities of the fat points in Z.
-- The conjecture is known to be correct for n<=9.

findhilb = (l) -> (
ww:=l;
if #l<3 then ww=join(l,{0,0,0});
n:=#ww;
n=n-1;
ww=zr(ww);
a1:=findalpha(ww);
tau:=findtau(ww);
d1:=a1-1;
<< "The output gives dim I_t, computed in degrees t from alpha(I)-1 to " << endl; 
<< "reg(I), where tau = reg(I)-1 is least degree such that" << endl;
<< "hilbert function of I equals hilbert polynomial of I." << endl;
if n>9 then (
   << "When more than 9 multiplicities are input," << endl;
   << "the output is a lower bound for dim I_t, which by the" << endl;
   << "SHGH conjecture should equal dim I_t." << endl);
<< endl;
<< " t      dim I_t" << "  (tau = " << tau << ")" << endl;
while (d1 <= tau+1) do (
   << " " << d1 << "     " << homcompdim(fundom({d1,ww})) << endl;
   d1=d1+1))

--input: l={n,m}, n = number of points, m = uniform multiplicity
--output: various bounds on alpha and tau

unifbounds = (l) -> (
n:=l#0;
m:=l#1;
ba:=bestrda(n);
bb:=bestrdb(n);
ea:=uniffindalpha(l);
t:=0;
<< "number of general points n of P2: " << n  << endl;
<< "multiplicity m of each point: " << m << endl;
<< endl;
if n<= 9 then (<< "Value of alpha: " << ea << endl)
else (
   << "Expected value of alpha (via SHGH conjecture): " << ea << endl;
   << " Note: The SHGH conjectural value of alpha is an upper bound." << endl);
<< "Lower Bounds on alpha:" << endl;
<< "  Roe's, via unloading: " << unifroealpha(l) << endl;
tmp:=unifezbhalpha(l);
<< "  Harbourne's, via Cor IV.i.2(a, b), using r="<<tmp#1<<" and d=";
<<tmp#2<<": "<< tmp#0 << endl;
t=ea;
while(t==unifbhalpha(l,ba#0,ba#1,t)) do t=t-1;
<< "  Harbourne's, via unloading, using r="<<ba#0<<" and d=";
<<ba#1<<": "<< t+1 << endl;
tmp=ezunifHRalpha(l);
<< "  Harbourne/Roe's first formula, using r="<<tmp#1<<" and d="<<
tmp#2<<": "<< tmp#0 << endl;
<< "  Harbourne/Roe's second formula, using r="<<(tmp#2)*(tmp#2)<<" and d="<<
tmp#2<<": "<<ezunifHRalphaB(l,(tmp#2)*(tmp#2),tmp#2) << endl;
r:=ba#0;
d:=ba#1;
t=ea;
while(t==unifHRalpha(l,r,d,t)) do t=t-1;
t=t+1;
tmp=ea;
while(tmp==unifHRalpha(l,bb#0,bb#1,tmp)) do tmp=tmp-1;
tmp=tmp+1;
if tmp>t then (
  t=tmp;
  r=bb#0;
  d=bb#1);
<<"  Harbourne/Roe's (via modified unloading), using r="<<r<<" and d="<<d;
<< ": "<<t<<endl<<endl;
tt:=uniffindtau(l);
if n<= 9 then ( << "Value of tau: " << tt << endl)
else (
   << "Expected value of tau (via SHGH conjecture): " << tt << endl;
   << " Note: The SHGH conjectural value of tau is a lower bound." << endl);
<< "Upper Bounds on tau:" << endl;
<< "  Hirschowitz's: " << Hiuniftau(l) << endl;
<< "  Gimigliano's: " << Guniftau(l) << endl;
if n>4 then ( << "  Catalisano's: " << Cuniftau(l) << endl);
t=0;
while(t*(t+3)-n*m*(m+1) < 2*t*(m-1)-2) do t=t+1;
<< "  Ballico's: " << t << endl;
t=0;
while(9*(t+3)*(t+3) <= 10*n*(m+1)*(m+1)) do t=t+1;
<< "  Xu's: " << t << endl;
t=0;
tmp=0;
while(tmp*tmp < n) do tmp=tmp+1;
while(2*t < 2*m*tmp + tmp - 3) do t=t+1;
<< "  Harbourne/Holay/Fitchett's: " << t << endl;
<< "  Roe's, via unloading: " << unifroetau(l) << endl;
tmp=ezunifHRtau(l);
<< "  Harbourne/Roe's first formula, using r="<<tmp#1<<" and d="<<
tmp#2<<": "<< tmp#0 << endl;
<< "  Harbourne/Roe's second formula, using r="<<(tmp#2)*(tmp#2)<<" and d="<<
tmp#2<<": "<<ezunifHRtauB(l,(tmp#2)*(tmp#2),tmp#2) << endl;
r=ba#0;
d=ba#1;
t=unifHRtau(l,r,d,tt);
tmp=unifHRtau(l,bb#0,bb#1,tt);
if tmp<t then (
  t=tmp;
  r=bb#0;
  d=bb#1);
<<"  Harbourne/Roe's (via unloading), using r="<<r<<" and d="<<d<< ": "<<t<<endl;
if (ba#0)*(bb#0)<(ba#1)*(bb#1)*n then (
   t = -3 + ceiling((m+1)*(bb#1)*n/(bb#0));
   r=bb#0;
   d=bb#1) else (
   t = -3 + ceiling((m+1)*(ba#0)/(ba#1));
   r=ba#0;
   d=ba#1);
<< "  Via Ran's observation, and Harbourne's bound on alpha," << endl;
<< "      using r="<<r<<" and d="<<d<< ": "<<t<<endl)

-- input: l={m1,...,mn}, n >=1 (number of points), m1, ... >=1 (the multiplicities)
-- output: various bounds on alpha and tau

bounds = (l) -> (
n:=#l;
ba:=bestrda(n);
bb:=bestrdb(n);
ea:=findalpha(zr(l));
t:=0;
tmp:=0;
<< "number of general points n of P2: " << n << endl;
<< "multiplicities of the points: " << l << endl;
<< endl;
if n<= 9 then ( << "Value of alpha: " << ea << endl)
else (
   << "Expected value of alpha (via SHGH conjecture): " << ea << endl;
   << " Note: The SHGH conjectural value of alpha is an upper bound." << endl);
<< "Lower bounds on alpha:" << endl;
<< "  Via Checking Psi: " << Psibound(l) << endl;
<< "  Roe's, via unloading: " << roealpha(l) << endl;
w:=ezbhalphaA(l);
<< "  Harbourne's, via Cor IV.i.2(a), using r="<<w#1<<", d=";
<<w#2<<": "<<w#0<<endl;
w=ezbhalphaB(l);
<< "  Harbourne's, via Cor IV.i.2(b), using r="<<w#1<<", d=";
<<w#2<<": "<<w#0<<endl;
w=ezbhalphaD(l);
<< "  Harbourne's, via Cor IV.i.2(d), using r="<<w#1<<", d="<<w#2;
<<", and j="<<w#3<<": "<<w#0<<endl;
r:=ba#0;
d:=ba#1;
t=ea;
while(t==HRalpha(l,r,d,t)) do t=t-1;
t=t+1;
tmp=ea;
while(tmp==HRalpha(l,bb#0,bb#1,tmp)) do tmp=tmp-1;
tmp=tmp+1;
if tmp>t then (
  t=tmp;
  r=bb#0;
  d=bb#1);
<<"  Harbourne/Roe's (via modified unloading), using r="<<r<<" and d="<<d;
<< ": "<<t<< endl << endl;
tt:=findtau(l);
if n<= 9 then ( << "Value of tau: " << tt << endl)
else (
   << "Expected value of tau (via SHGH conjecture): " << tt << endl;
   << " Note: The SHGH conjectural value of tau is a lower bound." << endl);
<< "Upper bounds on tau:" << endl;
<< "  Hirschowitz's: " << Hitau(l) << endl;
<< "  Gimigliano's: " << Gtau(l) << endl;
nn:=0;
i:=0;
w=prmt(zr(l));
scan(#l, i->(if w#i >0 then nn=nn+1));
if nn>4 then ( << "  Catalisano's: " << Ctau(l) << endl);
<< "  Roe's, via unloading: " << roetau(l) << endl;
r=ba#0;
d=ba#1;
t=HRtau(l,r,d,tt);
tmp=HRtau(l,bb#0,bb#1,tt);
if tmp<t then (
  t=tmp;
  r=bb#0;
  d=bb#1);
<< "  Harbourne/Roe's (via unloading), using r="<<r<<" and d="<<d<<": "<<t<<endl<<endl;
<<"One more bound on alpha; this one can be slow since it tries all r and d." << endl;
w=bestbhalpha(l,ea);
<< "  Harbourne's alpha lower bound (via unloading), using r=" << w#1 << " and d="; 
<< w#2 << ": " << w#0 << endl)

--decomp: prints a decomposition F=H+N for any divisor class F in Psi
--as described in the print statements of the script.
--input: decomp(l), where l={d,{m1,...,mn}}, signifying the divisor
--class F = dE_0-(m_1E_1+...+m_nE_n).

decomp = (l) -> (
<< "Let Psi be the subsemigroup of divisor classes generated by " << endl;
<< "exceptional classes and by -K. For any divisor class F, this " << endl;
<< "script determines if F is in Psi, and if so gives a decomposition" << endl;
<< "F=H+N, where N is a sum of exceptionals orthogonal to each other and to H" << endl;
<< "and H is in Psi but H.E >= 0 for all exceptionals E. The point of this" << endl;
<< "is that dim |F| = dim |H|, and conjecturally dim |H| = (H^2-H.K)/2." << endl << endl;
i:=0;
j:=0;
w:=l;
v:={};
ww:={};
ex:={};
mult:=0;
tmp:=fundom(l);
if tmp#0<tmp#1#0 or tmp#0<0 then (<< "Your class is not in Psi." << endl) else (
   if #l#1<3 then w={l#0,join(l#1,{0,0,0})};
   d:=3*(#(w#1));  --define an element {d,v} in fundamental domain
   scan(#(w#1), i->(v=join(v,{(#(w#1)) - i})));
   ww=fundomboth(w,{d,v});
   if ww#0#1 == zr(ww#0#1) then (<< "N = 0" << endl ) else (
   << "N is a sum of the following fixed exceptional classes:" << endl;
   if (ww#0#0) - (ww#0#1#0) - (ww#0#1#1) < 0 then (
       scan(#(w#1), j->(if j<=1 then ex=join(ex,{1}) else ex=join(ex,{0})));
       ex={1,ex};
       ex=(fundomboth(ww#1,ex))#1;
       <<ex<<" is a fixed component of multiplicity ";
       <<(ww#0#1#0)+(ww#0#1#1)-(ww#0#0)<<endl;
       ex={(ww#0#0)-(ww#0#1#1),(ww#0#0)-(ww#0#1#0)};
       scan(#(w#1), j->(if j>1 then ex=join(ex,{ww#0#1#j})));
       ww={{2*(ww#0#0)-(ww#0#1#0)-(ww#0#1#1),ex},ww#1});
   scan(#(w#1), i->(
         if (ww#0#1)#i<0 then (
            ex={};
            mult=-(ww#0#1)#i;
            scan(#(w#1), j->(if j==i then ex=join(ex,{-1}) else ex=join(ex,{0})));
            ex={0,ex};
            ex=(fundomboth(ww#1,ex))#1;
            << ex << " is a fixed component of multiplicity " << mult << endl))));
   << endl << "and H = " << (fundomboth(ww#1,{ww#0#0,zr(ww#0#1)}))#1 << endl))

-- input: list l of multiplicities for fat points Z
-- output: least t such that F_t(Z) is in Psi

Psibound = (l) -> (
t:=findalpha(l);
tmp:=fundom({t,l});
while(tmp#0>=tmp#1#0 and tmp#0>=0) do (t=t-1;
   tmp=fundom({t,l}));
t+1) 

-- prmt: arranges the elements of the list l={m_1,...,m_n} in descending order
-- Call it as prmt(l) where l is a list of integers.

prmt = (l) -> (  
(prmtboth(l,l))#0)

-- prmtboth: arranges the elements of the list l1 in descending order,
-- and applies the same permutation to l2
-- Call it as prmt(l1,l2) where l1 and l2 are lists of integers.

prmtboth = (l1,l2) -> (  
tmpv1:=l1;
tmpv2:=l2;
v1:=l1;
v2:=l2;
i:=0;
j:=0;
k:=0;
scan(#l1, i->(scan(#l1, j->(
            if tmpv1#i < tmpv1#j then (
               if i < j then (
                  k=-1;
                  v1={};
                  v2={};
                  while(k<#l1-1) do (k=k+1;
                     if k==i then (v1=join(v1,{tmpv1#j});
                        v2=join(v2,{tmpv2#j})) else (
                        if k==j then (v1=join(v1,{tmpv1#i});
                           v2=join(v2,{tmpv2#i})) else (v1=join(v1,{tmpv1#k});
                              v2=join(v2,{tmpv2#k}))));
                  tmpv1=v1;
                  tmpv2=v2))))));
{v1,v2})

-- zr: replaces negative values in a list l by zeroes.
-- Call it as zr(l) where l is a list of integers.

zr = (l) -> (    
v:={};
i:=0;
scan(#l, i->(
      if l#i<0 then v=join(v,{0}) else v=join(v,{l#i})));
v)

-- quad: performs a quadratic transform on a divisor class dE_0-(m_1E_1+...+m_nE_n).
-- Call it as quad({d,{m1,...,mn}}). The output is
-- {2d-m1-m2-m3,{d-m2-m3,d-m1-m3,d-m1-m2,m4,...,mn}}.

quad = (l) -> (
i:=0;
w:=l; 
if #l#1<3 then w={l#0,join(l#1,{0,0,0})};
v:={w#0 - w#1#1 - w#1#2,w#0 - w#1#0 - w#1#2,w#0 - w#1#0 - w#1#1};
scan(#w#1, i->(if i>2 then v=join(v,{w#1#i})));
v={2*(w#0) - w#1#0 - w#1#1 - w#1#2,v};
v)

-- fundom: Call it as fundom({d,{m1,...,mn}}). The output is a new 
-- list {d',{m1',...,mn'}}; the class dE_0-(m_1E_1+...+m_nE_n) is
-- equivalent via Cremona transformations to d'E_0-(m_1'E_1+...+m_n'E_n),
-- where d' is either negative or as small as possible.

fundom = (l) -> (
(fundomboth(l,l))#0)

-- fundomboth: applies fundom to l1 to reduce l1 to fundamental
-- domain of a certain group operation, and applies the same
-- group operation g to l2. If l2 starts out in the fundamental domain,
-- and {l1',l2'}=fundomboth(l1,l2), then {l2,g^{-1}l}=fundomboth(l2',l).
-- This allows one to compute the action of g^{-1}.

fundomboth = (l1,l2) -> (
w1:=l1; 
w2:=l2;
v:={}; 
if #l1#1<3 then w1={l1#0,join(l1#1,{0,0,0})};
if #l2#1<3 then w2={l2#0,join(l2#1,{0,0,0})};
v=prmtboth(w1#1,w2#1);
w1={w1#0,v#0};
w2={w2#0,v#1};
while ((w1#0 < w1#1#0 + w1#1#1 + w1#1#2) and (w1#0 >= 0)) do (
   w1=quad(w1);
   w2=quad(w2);
   v=prmtboth(w1#1,w2#1);
   w1={w1#0,v#0};
   w2={w2#0,v#1});
{w1,w2})

-- homcompdim: computes e(F_t(Z)), the expected dimension of a component I_d 
-- of a fat points ideal I corresponding to a fat point subscheme Z of general
-- points taken with multiplicities m_1, ..., m_n. Call it as 
-- homcompdim({d,{m_1,...,m_n}}); the output is the SHGH conjectural 
-- dimension of I_d, which is the actual dimension if n < 10.

homcompdim = (l) -> (
h:=0;
i:=0;
w:=l; 
if #l#1<3 then w={l#0,join(l#1,{0,0,0})};
w=fundom(w);
d:=w#0;
w=fundom({d,zr(w#1)});
d=w#0;
v:=zr(w#1);
if d<0 then h=0 else (
   tmp:=0;
   scan(#v, i->(tmp = v#i*v#i+v#i+tmp));
   h=floor((d*d+3*d+2-tmp)/2);
   if h < 0 then h=0);
h)

-- findalpha: find alpha, the least degree t such that
-- I_t \ne 0, where I is the ideal corresponding to n general
-- points taken with multiplicities m_1, ...,m_n. Call it as 
-- findalpha({m_1,...,m_n}). The output is the SHGH conjectural value 
-- of alpha, which is the actual value if n < 10 and an upper bound otherwise.

findalpha = (l) -> (
i:=1;
w:=prmt(zr(l)); 
if #l<3 then w=join(l,{0,0,0});
d:=w#0;                -- alpha is at least the max mult
if (#w)<9 then (       -- if n<=8, to speed things, make an estimate
  while(i < (#w)) do (d=d+w#i;
     i=i+1);
  d=ceiling(d/3);
  if d < w#0 then d = w#0);
while (homcompdim({d,w}) < 1) do d=d+1;
d)

-- dot: computes a dot product of two lists l1 and l2 (of equal length)
-- of integers. Call it as dot(l1,l2).

dot = (l1,l2) -> (
i:=0;
dottot:=0;
scan(#l1, i->(dottot = dottot + (l1#i)*(l2#i)));
dottot)

-- input: l={m1,...,mn}, n >=1 (number of points), m1,... (the multiplicities)
-- output: the SHGH conjectured value of tau; this is the actual value if 
-- n < 10, and a lower bound otherwise.

findtau = (l) -> (
t:=findalpha(l);
if t > 0 then t = t-1;  -- tau is at least alpha - 1
n:=#l;
v:=l;
p:=dot(v,v);
K:={};
j:=0;
q:=0;
scan(#l, j->(q=q+l#j));
while(2*(homcompdim(join({t},{v}))) > t*t-p+3*t-q+2) do t=t+1;
t)

-- input: positive integer n
-- output: {r,d}, where r^2>=d^2n, n>=r, and nd/r is as big as possible

bestrda = (n) -> (
rootn:=0;
while(rootn*rootn<=n) do rootn=rootn+1;
rootn=rootn-1;
d:=1;
r:=0;
if rootn*rootn==n then r=rootn else r=rootn+1;
tmpr:=1;
tmpd:=1;
while(tmpd<=rootn) do (
   tmpr=tmpd*rootn;
   while(tmpr*tmpr<tmpd*tmpd*n) do tmpr=tmpr+1;
   if tmpr*d < tmpd*r then (
      r=tmpr;
      d=tmpd);
   tmpd=tmpd+1);
{r,d})

-- input: positive integer n
-- output: {r,d}, where r^2<=d^2n, n>=r, and r/d is as big as possible

bestrdb = (n) -> (
rootn:=0;
while(rootn*rootn<n) do rootn=rootn+1;
d:=1;
r:=0;
if rootn*rootn==n then r=rootn else r=rootn-1;
tmpr:=1;
tmpd:=1;
while(tmpd<=rootn) do (
   tmpr=tmpd*rootn;
   while(tmpr*tmpr>tmpd*tmpd*n) do tmpr=tmpr-1;
   if tmpr>n then tmpr=n;
   if tmpr*d > tmpd*r then (
      r=tmpr;
      d=tmpd);
   tmpd=tmpd+1);
{r,d})

-- input: l={m1,...,mn}, n= number of points, mi = multiplicity of ith point
-- output: Roe's algorithmic lower bound on alpha

roealpha = (l) -> (
i:=0;
v={};
w2:={};
w:=l;
i1:=2;
if #l<3 then w=join(w,{0,0,0});
while(i1<#w) do (
   v={1};
   scan(#w, i->(if i>0 then (if i<=i1 then v=join(v,{-1}) else v=join(v,{0}))));
   w=zr(prmt(w));
   while(dot(w,v)<0) do (
      w2={};
      scan(#w, i->(w2 = join(w2,{w#i+v#i})));
      w=zr(prmt(w2)));
   i1=i1+1);
w#0)

-- input: l={n,m}, n >=1 (number of points), m >=1 (uniform multiplicity)
-- output: Roe's algorithmic lower bound on alpha

unifroealpha = (l) -> (
i1:=0;
intchk:=0;
n:=l#0;
m:=l#1;
roebnd:=m;
q:=n-1; -- q keeps track during unloading of number of 
        -- points after the first with maximum multiplicity
if n>2 then (
   i1=2;
   while(i1<n) do (
      if q<i1 then intchk=i1*m-(i1-q) else intchk=i1*m;
      while(roebnd<intchk) do (
         roebnd=roebnd+1;
         if q<i1 then (q=n-i1+q-1;
            m=m-1) else (q=q-i1;
            if q==0 then (m=m-1;
               q=n-1));
         if q<i1 then intchk=i1*m-(i1-q) else intchk=i1*m);
      i1=i1+1)) else roebnd = m;
roebnd)

-- input: l={m1,...,mn}, n= number of points, mi = multiplicity of ith point
-- output: Roe's algorithmic upper bound on tau

roetau = (l) -> (
i:=0;
vv={};
vv1:={};
w2:={};
ww:=l;
i1:=1;
while(i1<#ww-1) do (
   vv={1};
   scan(#ww, i->(if i>0 then (if i-1<=i1 then vv=join(vv,{-1}) else vv=join(vv,{0}))));
   vv1={1};
   scan(#ww, i->(if i>0 then (if i<=i1 then vv1=join(vv1,{-1}) else vv1=join(vv1,{0}))));
   ww=zr(prmt(ww));
   while((dot(ww,vv)) < -1) do (
      w2={};
      scan(#ww, i->(w2 = join(w2,{ww#i+vv1#i})));
      ww=zr(prmt(w2)));
   i1=i1+1);
ww#0+ww#1-1)

-- input: l={n,m}, n >=1 (number of points), m >=1 (uniform multiplicity)
-- output: Roe's algorithmic upper bound on tau

unifroetau = (l) -> (
i:=1;
s:=0;
n:=l#0;
m1:=l#1;
m2:=0;
if n>1 then m2=m1;
n2:=n-1; -- n2 keeps track during unloading of number of 
         -- points with multiplicity m2
while(i < (n-1)) do (
   if (i+1) <= n2 then s=(i+1)*m2 else s=(i+1)*m2+n2-i-1;
   while((m1-s) < -1) do (
      m1=m1+1;
      if i < n2 then n2=n2-i else (
         n2=n-i+n2-1;
         m2=m2-1);
      if (i+1) <= n2 then s=(i+1)*m2 else s=(i+1)*m2+n2-i-1);
   i=i+1);
m1+m2-1)

-- input: l={n,m}, n >=1 (number of points), m >=1 (uniform multiplicity)
-- output: list {i,r,d}, with i being Harbourne's easy lower bound on 
-- alpha (via Cor IV.i.2 (a), (b)) computed using the best r and d

unifezbhalpha = (l) -> (
w:=bestrda(l#0);
i:=ceiling((l#1)*(l#0)*(w#1)/(w#0));  -- compute mnd/r rounded up
r:=w#0;
d:=w#1;
w=bestrdb(l#0);
j:=ceiling((l#1)*(w#0)/(w#1));  -- compute mr/d rounded up
if i<j then (i=j;
  r=w#0;
  d=w#1);
{i,r,d})

-- input: {n,m}, n >=1 (number of points), m >=1 (uniform multiplicity)
-- output: {a,r,d}, where a is lower bound on alpha via formula of [HR] 
-- using r and d. 

ezunifHRalpha = (l) -> (
n:=l#0;
m:=l#1;
t:=0;
d:=0;
while (d*d <= n) do d=d+1;
d=d-1;
r:=d;
while (r*r < d*d*n) do r=r+1;
q:=ceiling(n*m/r)-1;
while(((t+2)*(t+1)<=2*(m*n-r*q))and t<d) do t=t+1;
{t+q*d,r,d})

-- input: ({n,m},r,d), n >=1 (number of points), m >=1 (uniform multiplicity)
-- where r <= n, and d(d+1)/2 <= r <= d^2
-- output: lower bound on alpha via formula of [HR]. 

ezunifHRalphaB = (l,r,d) -> (
n:=l#0;
m:=l#1;
g:=(d-1)*(d-2)/2;
tmp:=floor((m*r+g-1)/d);
t:=ceiling(n*m/r)-1;
rr:=m*n-r*t;
s:=0;
while(((s+1)*(s+2) <= 2*rr) and s<d) do s=s+1;
s=s-1;
t=s+t*d;
if tmp<t then t=tmp;
t+1)

-- input: ({n,m},r,d), n >=1 (number of points), m >=1 (uniform multiplicity)
-- where r <= n, and r <= d^2
-- output: upper bound on tau via formula of [HR]. 

ezunifHRtauB = (l,r,d) -> (
n:=l#0;
m:=l#1;
g:=(d-1)*(d-2)/2;
q:=ceiling(n*m/r)-1;
rr:=m*n-r*q;
t:=q*d+ceiling((rr+g-1)/d);
tmp:=q*d+d-2;
if t<tmp then t=tmp;
t)

-- input: ({n,m},r,d,ea), n >=1 (number of points), m >=1 (uniform multiplicity)
-- ea is an estimate for alpha (for speed); it must be set equal to
-- a value no bigger than the eventual value of unifHRalpha (0, for example)
-- output: lower bound on alpha via modified unloading method of [HR], 
-- computed using given r and d

unifHRalpha = (l,r,d,ea) -> (
i:=ea-1;
n2:=0;  -- n2 keeps track of the number of points with maximum multiplicity
s:=0;
tmpi:=-1;
tmpm:=0;
g:=(d-1)*(d-2)/2;
while(tmpi<tmpm) do (
   i=i+1;
   tmpi=i;
   tmpm=l#1;
   n2=l#0;
   if r <= n2 then s=r*tmpm else s=r*tmpm-r+n2;
   while(((tmpi*d-s < g) and (tmpi >= d-2)) or (((tmpi+1)*(tmpi+2)<=2*s) and 
                                               (tmpi<d) and (tmpi>=0))) do (
         tmpi=tmpi-d;       
         if r<n2 then n2=n2-r else (          
             n2=l#0-r+n2;
             tmpm=tmpm-1;
             if tmpm <= 0 then (
                tmpm=0;
                n2=l#0));
         if r<=n2 then s=r*tmpm else s=r*tmpm-r+n2));
i)

-- input: ({m_1,...,m_n},r,d,ea), n >=1 (number of points), m_i >=1 (multiplicities)
-- ea is an estimate for alpha (for speed); it must be set equal to
-- a value no bigger than the eventual value of HRalpha (0, for example)
-- output: lower bound on alpha via modified unloading method of [HR], 
-- computed using given r and d

HRalpha = (l,r,d,ea) -> (
i:=ea-1;
j:=0;
n:=0;
g:=(d-1)*(d-2)/2;
v:=prmt(zr(l));
tmpv:=v;
ttmpv:={};
scan(#l, j->(if v#j>0 then n=n+1));
if n>0 then (
   ww:={};
   scan(#l, j->(if j<r then ww=join(ww,{1}) else ww=join(ww,{0})));
   tmpi:=-1;
   while(tmpi < tmpv#0) do (
      i=i+1;
      tmpi=i;
      tmpv=v;
      while(((tmpi*d-(dot(ww,tmpv))<g) and (tmpi >= d-2)) or 
            (((tmpi+1)*(tmpi+2)<=2*(dot(ww,tmpv))) and (tmpi<d) and (tmpi>=0))) do (
         tmpi=tmpi-d;
         ttmpv={};
         scan(#l, j->(ttmpv=join(ttmpv,{(tmpv#j)-(ww#j)})));
         tmpv=prmt(zr(ttmpv)))));
i)

-- input: ({n,m},r,d,ea), n >=1 (number of points), m >=1 (uniform multiplicity)
-- ea is an estimate for alpha (for speed); it must be set equal to
-- a value no bigger than the eventual value of unifbhalpha (0, for example)
-- output: Harbourne's algorithmic lower bound on alpha via unloading, 
-- using the given r and d. 

unifbhalpha = (l,r,d,ea) -> (
i:=ea-1;
n2:=0;  -- n2 keeps track of the number of points with maximum multiplicity
s:=0;
tmpi:=-1;
tmpm:=0;
while(tmpi<0) do (
   i=i+1;
   tmpi=i;
   tmpm=l#1;
   n2=l#0;
   if r <= n2 then s=r*tmpm else s=r*tmpm-r+n2;
   while((tmpi*d-s < 0) and (tmpi >= 0)) do (
      tmpi=tmpi-d;
      if r<n2 then n2=n2-r else (
         n2=l#0-r+n2;
         tmpm=tmpm-1;
         if tmpm <= 0 then (
            tmpm=0;
            n2=l#0));
      if r<=n2 then s=r*tmpm else s=r*tmpm-r+n2));
i)

-- input: ({n,m},r,d), n >=1 (number of points), m >=1 (uniform multiplicity)
-- where r <= n, and 2r >= n+d^2
-- output: lower bound on alpha via formula which agrees with that via unloading. 

ezunifBHalphaB = (l,r,d) -> (
n:=l#0;
m:=l#1;
q:=floor(n*m/r);
rr:=m*n-r*q;
t:=q-1+ceiling(rr/d);
tmp:=d*ceiling(n*m/r);
if tmp<t then t=tmp;
t+1)

-- input: l={m1,...,mn}, n >=1 (number of points), m1, ... >=1 (the multiplicities)
-- output: list (aa,r,d), with aa being Harbourne's lower bound on 
-- alpha (via Cor IV.i.2(a)), computed using the best r and d

ezbhalphaA = (l) -> (
i:=0;
s:=0;
n:=0;
v:=prmt(zr(l));
scan(#l, i->(if v#i>0 then n=n+1));
w:=bestrda(n);
i=0;
while(i<n) do (
  s=s+v#i;
  i=i+1);
best:=ceiling(s*(w#1)/(w#0));
{best, w#0, w#1})

-- input: l={m1,...,mn}, n >=1 (number of points), m1, ... >=1 (the multiplicities)
-- output: list (aa,r,d), with aa being Harbourne's lower bound on 
-- alpha (via Cor IV.i.2(b)), computed using the best r and d

ezbhalphaB = (l) -> (
i:=0;
s:=0;
n:=0;
v:=prmt(zr(l));
scan(#l, i->(if v#i>0 then n=n+1));
w:=bestrdb(n);
i=0;
while(i<n) do (
  s=s+v#i;
  i=i+1);
best:=ceiling(s*(w#0)/(n*(w#1)));
{best, w#0, w#1})

-- input: l={m1,...,mn}, n >=1 (number of points), m1, ... >=1 (the multiplicities)
-- output: list (aa,r,d,j), with aa being Harbourne's lower bound on 
-- alpha (via Cor IV.i.2(d)), computed using the best r, d and j

ezbhalphaD = (l) -> (
w:={};
i:=0;
scan(#l, i->(if l#i>0 then w=join(w,{l#i})));
w=prmt(w);
n:=#w;
bnd:=0;
tmpbnd:=0;
r:=0;
d:=0;
j:=0;
tmpr:=0;
tmpd:=0;
tmpj:=0;
if n>0 then (
   while(tmpr<n) do (tmpr=tmpr+1;
      tmpd=0;
      while(tmpd*tmpd<tmpr) do (tmpd=tmpd+1;
         tmpj=0;
         while(tmpj<tmpd*tmpd) do (tmpj=tmpj+1;
             tmpbnd=lpa(w,tmpr,tmpd,tmpj);
             if tmpbnd>bnd then (
                 bnd=tmpbnd;
                 r=tmpr;
                 d=tmpd;
                 j=tmpj)))));
{bnd,r,d,j})

-- lpa computes the bound given in Cor IV.i.2(d); attempts
-- various solutions with the hope of approximating the optimal 
-- solution to the linear programming problem indicated by Thm IV.i.1.
-- lpa is called by ezbhalphaD

lpa = (l,r,d,j) -> (
i:=0;
n:=#l;
sum:=0;
sumb:=0;
bnd:=0;
if d*d >= r then  (
   scan(#l, i->(if i<r then sum=sum+l#i));
   bnd=ceiling(sum/d)) else (
   if j==0 then (
      scan(#l, i->(if i<d*d then sum=sum+l#i));
      bnd=ceiling(sum/d)) else (
      M:=floor((r-d*d)*(r-d*d+j)/j);
      scan(#l, i->(if i<d*d-j then sum=sum+l#i else (if i<M+r then sumb=sumb+l#i)));
      sumb=sumb*j/(r-d*d+j);
      if M<n-r then sumb=sumb+(l#(M+r))*(r-d*d-j*M/(r-d*d+j));
      bnd=ceiling((sumb+sum)/d)));
bnd)

-- input: ({m1,...,mn},r,d,ea), n >=1 (number of points), m1, ... >=1 
-- (the multiplicities), r and d positive integers, ea any value
-- not bigger than the eventual value of bhalpha; can be set to 0 
-- output: Harbourne's unloading lower bound on alpha, using given r and d

bhalpha = (l,r,d,ea) -> (
i:=ea-1;
j:=0;
n:=0;
v:=prmt(zr(l));
tmpv:=v;
ttmpv:={};
scan(#l, j->(if v#j>0 then n=n+1));
if n>0 then (
   ww:={};
   scan(#l, j->(if j<r then ww=join(ww,{1}) else ww=join(ww,{0})));
   tmpi:=-1;
   while(tmpi < tmpv#0) do (
      i=i+1;
      tmpi=i;
      tmpv=v;
      while((tmpi*d-(dot(ww,tmpv))<0) and (tmpi >= tmpv#0)) do (
         tmpi=tmpi-d;
         ttmpv={};
         scan(#l, j->(ttmpv=join(ttmpv,{(tmpv#j)-(ww#j)})));
         tmpv=prmt(zr(ttmpv)))));
i)

-- Find bhalpha using best possible r and d;
-- ea is an a priori estimate for alpha (for speed)
-- it must be set to a value >= than the actual value 
-- of alpha (e.g., ea=findalpha(l))

bestbhalpha = (l,ea) -> (
w:={};
i:=0;
scan(#l, i->(if l#i>0 then w=join(w,{l#i})));
w=prmt(w);
n:=#w;
bnd:=0;
tmpbnd:=0;
r:=0;
d:=0;
tmpr:=0;
tmpd:=0;
if n>0 then (
   while(tmpr<n) do (tmpr=tmpr+1;
      tmpd=0;
      while(tmpd*tmpd<tmpr) do (tmpd=tmpd+1;
             tmpbnd=ea;
             while(tmpbnd==bhalpha(w,tmpr,tmpd,tmpbnd)) do tmpbnd=tmpbnd-1;
             tmpbnd=tmpbnd+1;
             if tmpbnd>bnd then (
                 bnd=tmpbnd;
                 r=tmpr;
                 d=tmpd))));
{bnd,r,d})

-- input: ({m1,...,mn},r,d,et), n >=1 (number of points), m1, ... >=1 (the multiplicities),
-- r and d positive integers, et any lower bound for tau (used for speed; can be set to 0).
-- output: Harbourne/Roe's algorithmic upper bound on tau, 
-- with given r and d (assumes char = 0).

HRtau = (l,r,d,et) -> (
i:=et-1;
j:=0;
n:=0;
v:=prmt(zr(l));
tmpv:=v;
ttmpv:={};
scan(#l, j->(if v#j>0 then n=n+1));
if n>0 then (
   ww:={};
   g:=(d-1)*(d-2)/2; -- genus of plane curve of degree d
   scan(#l, j->(if j<r then ww=join(ww,{1}) else ww=join(ww,{0})));
   tmpi:=0;
   while(tmpv#0 > 0) do (
      i=i+1;
      tmpi=i;
      tmpv=v;
      while((tmpi*d-(dot(ww,tmpv))>=g-1) and (tmpi>=d-2) and (tmpv#0 >0)) do (
         tmpi=tmpi-d;
         ttmpv={};
         scan(#l, j->(ttmpv=join(ttmpv,{(tmpv#j)-(ww#j)})));
         tmpv=prmt(zr(ttmpv)))));
i)

-- input: ({n,m},r,d,et), n >=1 (number of points), m >=1 (the uniform multiplicity),
-- r and d positive integers, et any lower bound for tau (used for speed; can be set to 0)
-- output: Harbourne/Roe's algorithmic upper bound on tau, using given r and d
-- (assumes char = 0).

unifHRtau = (l,r,d,et) -> (
i:=et-1;
n2:=0;  -- n2 is the number of points with maximum multiplicity
s:=0;
tmpm:=1;
tmpi:=0;
g:=(d-1)*(d-2)/2; -- genus of plane curve of degree d
while(tmpm > 0) do (
   i=i+1;
   tmpi=i;
   tmpm=l#1;
   n2= l#0;
   if r<=n2 then s=r*tmpm else s=r*tmpm-r+n2;
   while((tmpi*d-s>=g-1) and (tmpi>=d-2) and (tmpm >0)) do (
      tmpi=tmpi-d;
      if r<n2 then n2=n2-r else (
          n2=l#0-r+n2;
          tmpm=tmpm-1;
          if tmpm<0 then (
              tmpm=0;
              n2= l#0));
          if r<=n2 then s=r*tmpm else s=r*tmpm-r+n2));
i)

-- input: l={n,m}, n >=1 (number of points), m >=1 (the uniform multiplicity)
-- output: list {a,r,d}, where a is Harbourne/Roe's formulaic upper 
-- bound on tau (char 0) computed using r and d

ezunifHRtau = (l) -> (
n:=l#0;
m:=l#1;
d:=0;
while (d*d <= n) do d=d+1;
d=d-1;
r:=d;
while (r*r < d*d*n) do r=r+1;
g:=(d-2)*(d-1)/2;
a:= ceiling((m*r+g-1)/d);
b:=-2+d*ceiling(m*n/r);
if a<b then a=b;
{a,r,d})

-- input: l={n,m}, n >=1 (number of points), m >=1 (the uniform multiplicity)
-- output: the SHGH conjectured value of alpha; this is the actual value if 
-- n < 10, and an upper bound otherwise.

uniffindalpha = (l) -> (
n:=l#0;
m:=l#1;
a:=-1;
if n==1 then a=m;
if n==2 then a=m;
if n==3 then a=ceiling(3*m/2);
if n==4 then a=2*m;
if n==5 then a=2*m;
if n==6 then a=ceiling(12*m/5);
if n==7 then a=ceiling(21*m/8);
if n==8 then a=ceiling(48*m/17);
if n==9 then a=3*m;
if n>9 then (
  while(a*a-n*m*m+3*a-n*m+2 <0) do a=a+m;
  a=a-m;
  while(a*a-n*m*m+3*a-n*m+2 <=0) do a=a+1);
a)

-- input: l={n,m}, n >=1 (number of points), m >=1 (the uniform multiplicity)
-- output: the SHGH conjectured value of tau; this is the actual value if 
-- n < 10, and a lower bound otherwise.

uniffindtau = (l) -> (
n:=l#0;
m:=l#1;
t:=-1;
if n==1 then t=m-1;
if n==2 then t=2*m-1;
if n==3 then t=2*m-1;
if n==4 then t=2*m;
if n==5 then t=ceiling((5*m-1)/2);
if n==6 then t=ceiling((5*m-1)/2);
if n==7 then t=ceiling((8*m-1)/3);
if n==8 then t=ceiling((17*m-1)/6);
if n==9 then t=3*m;
if n>9 then (
  while(t*t-n*m*m+3*t-n*m+2 <0) do t=t+m;
  t=t-m;
  while(t*t-n*m*m+3*t-n*m+2 <0) do t=t+1);
if t<0 then t=0;
t)

-- input: l={n,m}, n >=1 (number of points), m (the multiplicity of each point)
-- output: Hirschowitz's lower bound for tau

Hiuniftau = (l) -> (
n:=l#0;
m:=l#1;
t:=m;
s:=n*m*(m+1);
a:=ceiling((t+3)/2);
b:=ceiling((t+2)/2);
while(a*b*2 <= s) do (
  t=t+1;
  a=ceiling((t+3)/2);
  b=ceiling((t+2)/2));
t)

-- input: l={m1,...,mn}, n >=1 (number of points), m1, ... >=1 (the multiplicities)
-- output: Hirschowitz's lower bound for tau

Hitau = (l) -> (
n:=#l;
i:=0;
w:=prmt(zr(l));
t:=w#0;
s:=0;
scan(#l, i->(s=s+(w#i)*((w#i)+1)));
a:=ceiling((t+3)/2);
b:=ceiling((t+2)/2);
while(a*b*2 <= s) do (
  t=t+1;
  a=ceiling((t+3)/2);
  b=ceiling((t+2)/2));
t)

-- input: l={n,m}, n >=1 (number of points), m (the multiplicity of each point)
-- output: Gimigliano's lower bound for tau

Guniftau = (l) -> (
n:=l#0;
m:=l#1;
t:=0;
while(t*(t+3)<2*n) do t=t+1;
m*t)

-- input: l={m1,...,mn}, n >=1 (number of points), m1, ... >=1 (the multiplicities)
-- output: Gimigliano's lower bound for tau

Gtau = (l) -> (
n:=0;
w:=prmt(zr(l));
scan(#l, i->(if w#i >0 then n=n+1));
t:=0;
s:=0;
i:=0;
while(t*(t+3)<2*n) do t=t+1;
scan(#l, i->(if i<t then s=s+w#i));
s)

-- input: l={n,m}, n >=5 (number of points), m>0 (the multiplicity of each point)
-- output: Catalisano's lower bound for tau

Cuniftau = (l) -> (
s:=l#0;
m:=l#1;
r:=0;
t:=0;
f:=0;
while(f*(f+1) <= 2*s) do f=f+1;
f=f-1;
while(2*r<2*s-f*(f+1)) do r=r+1;
d1:=0;
d:=f;
if r==0 then d1=f-1 else d1=f;
t=d1+(m-1)*d;
if 2*t+1 < 5*m then t=ceiling((5*m-1)/2);
if t<2*m-1 then t=2*m-1;
if r == f then (if s >= 9 then t=m*d1+1);
t)

-- input: l={m1,...,mn}, n >=5 (number of points), m1, ... >=1 (the multiplicities)
-- output: Catalisano's lower bound for tau

Ctau = (l) -> (
n:=0;
i:=0;
w:=prmt(zr(l));
scan(#l, i->(if w#i >0 then n=n+1));
vm:={};
vs:={};
i=0;
while(i < n-1) do (
  if w#i > w#(i+1) then (
  vs=join(vs,{i+1});
  vm=join(vm,{w#i}));
  i=i+1);
vs=join(vs,{n});
vm=join(vm,{w#(n-1)});
i=#vm - 1;
v:={vm#i};
while(i > 0) do (
  v=join({vm#(i-1) - vm#i},v);
  i=i-1);
vf:={};
vr:={};
scan(#vs, i->(
  f:=0;
  r:=0;
  while(f*(f+1) <= 2*(vs#i)) do f=f+1;
  f=f-1;
  while(2*r<2*(vs#i)-f*(f+1)) do r=r+1;
  vf=join(vf,{f});
  vr=join(vr,{r})));
t:=0;
if (vr#(#vr-1)) == 0 then t = - 1;
d1:=t+vf#(#vf-1);
scan(#vs, i->(
  t=t+(vf#i)*(v#i)));
if 2*t+1 < (w#0)+(w#1)+(w#2)+(w#3)+(w#4) then 
t=ceiling(((w#0)+(w#1)+(w#2)+(w#3)+(w#4)-1)/2);
if t<(w#0)+(w#1)-1 then t=(w#0)+(w#1)-1;
if (vr#0) == (vf#0) then (if s >= 9 then (if (w#0)==(w#(n-1)) 
    then (if (w#0)>1 then t=(w#0)*d1+1)));
if (vr#0) == 0 then (if s > 9 then (if (w#0)==(w#(n-2)) 
    then (if (w#(n-1)) == 1 then t=(w#0)*d1+1)));
t)

\endverbatim
\parindent=20pt

\tenpoint


\References

\bibitem{\refAH} Alexander, J. and Hirschowitz, A. {\it Polynomial interpolation 
in several variables}, 
J.\ Alg.\ Geom. 4 (1995), 201--222.

\bibitem{\refAHb} Alexander, J. and Hirschowitz, A. {\it An asymptotic 
vanishing theorem for generic unions of multiple points}, 
Invent. Math. 140 (2000), no. 2, 303--325.

\bibitem{\refBal} Ballico, E. {\it Curves of minimal degree
with prescribed singularities}, Illinois J.\ Math.\ 45 (1999), 672--676.

\bibitem{\refCam} Campanella, G. {\it Standard 
bases of perfect homogeneous
polynomial ideals of height $2$}, 
J.\ Alg.\  101 (1986), 47--60.

\bibitem{\refCatb} Catalisano, M.\ V. {\it Linear Systems of Plane Curves 
through Fixed ``Fat'' Points of \pr2}, 
J.\ Alg.\  142 (1991), 81-100.

\bibitem{\refCat} Catalisano, M.\ V. {\it ``Fat'' points on a conic}, 
Comm.\ Alg.\  19(8) (1991), 2153--2168.

\bibitem{\refCMa} Ciliberto, C. and Miranda, R. {\it Degenerations
of planar linear systems}, 
J. Reine Angew. Math. 501 (1998), 191-220.

\bibitem{\refCMb} Ciliberto, C. and Miranda, R. {\it Linear systems
of plane curves with base points of equal multiplicity}, 
Trans. Amer. Math. Soc. 352 (2000), 4037--4050.

\bibitem{\refDGM} Davis, E.\ D., Geramita, A.\ V., and 
Maroscia, P. {\it Perfect
Homogeneous Ideals: Dubreil's Theorems Revisited},
Bull.\ Sc.\ math., $2^e$ s\'erie, 108 (1984), 143--185.

\bibitem{\refDu} Dubreil, P.
{\it Sur quelques propri\'et\'es des syst\`emes de points dans
le plan et des courbes gauches alg\'ebriques},
Bull.\ Soc.\ Math.\ France, 61 (1933), 258--283.

\bibitem{\refDuVa} Du Val, P.
{\it The unloading problem for plane curves},
Amer. J. Math. 62(1940), 307--311.

\bibitem{\refDuVb} Du Val, P.
{\it Application des id\'ees cristallographiques 
\`a l'\'etude des groupes de transformations cr\'emoniennes},
1960 3 i\`eme Coll. G\'eom. Alg\'ebrique (Bruxelles, 1959) pp. 65--73 
Centre Belge Rech. Math., Louvain.

\bibitem{\refEC} Enriques, F. and Chisini, O.
{\it Lezioni sulla teoria geometrica delle
  equazioni e delle funzioni algebriche}.
N. Zanichelli, Bologna, 1915.

\bibitem{\refE} Evain, L.
{\it La fonction de Hilbert de la r\'eunion
de $4^h$ gros points g\'en\'eriques
de \pr2 de m\^eme multiplicit\'e},
J. Alg. Geom. 8 (1999), 787--796.

\bibitem{\refFL} Fatabbi, G. and Lorenzini, A.
{\it On a sharp bound for the regularity index
of any set of fat points},
preprint (2000).

\bibitem{\refFione} Fitchett, S. {\it On bounding the number
of generators for fat point ideals on the projective plane}, 
to appear, J. Alg.

\bibitem{\refFitwo} Fitchett, S. {\it Maps of linear systems
on blow ups of the projective plane}, 
to appear, J. Pure and Applied Alg.

\bibitem{\refFHH} Fitchett, S., Harbourne, B., and Holay, S.
{\it Resolutions of Fat Point Ideals involving 8 General Points of \pr2}, 
preprint, 2000.

\bibitem{\refGGR} Geramita, A.\ V., Gregory, D.\ and Roberts, L.
{\it Monomial ideals and points in projective space},
J.\ Pure and Appl.\ Alg. 40 (1986), 33--62.

\bibitem{\refGiThesis} Gimigliano, A.
{\it On Linear Systems of Plane Curves},
Ph. D. thesis, Queen's University, Kingston, Ontario (1987).

\bibitem{\refGi} Gimigliano, A. {\it Our thin knowledge of fat points},
Queen's papers in Pure and Applied Mathematics, no. 83,
The Curves Seminar at Queen's, vol. VI (1989).

\bibitem{\refGib} Gimigliano, A.
{\it Regularity of Linear Systems of Plane Curves},
J. Alg. 124 (1989), 447--460.

\bibitem{\refHHF} Harbourne, B., Holay, S. and Fitchett, S.
{\it Resolutions of Ideals of Quasiuniform Fat Point Subschemes of \pr2},
preprint (2000).

\bibitem{\refvanc} Harbourne, B. {\it The geometry of 
rational surfaces and Hilbert
functions of points in the plane},
Can.\ Math.\ Soc.\ Conf.\ Proc.\ 6 
(1986), 95--111.

\bibitem{\refsundance} \manyby. {\it Iterated blow-ups and moduli for 
rational surfaces}, in: Algebraic Geometry, 
Sundance 1986, LNM \#1311, (1988), 101--117.

\bibitem{\refravello} \manyby. {\it 
Points in Good Position in \pr 2}, in:
Zero-dimensional schemes, Proceedings of the
International Conference held in Ravello, Italy, June 8--13, 1992,
De Gruyter, 1994.

\bibitem{\refpos} \manyby. {\it Rational Surfaces with $K^2>0$}, 
Proc. Amer. Math. Soc. 124, 727--733 (1996).

\bibitem{\refanti} \manyby. {\it Anticanonical rational surfaces}, 
Trans. Amer. Math. Soc. 349, 1191--1208 (1997).

\bibitem{\reffatpts} \manyby. {\it Free Resolutions of Fat Point 
Ideals on \pr2}, J. Pure and Applied Alg. 125 
(1998), 213--234.

\bibitem{\refigp} \manyby. {\it  The Ideal Generation 
Problem for Fat Points},
J. Pure and Applied Alg. 145 (2000), 165--182.

\bibitem{\refseven} \manyby. {\it  An Algorithm for Fat Points on \pr2},
Canad. J. Math. 52 (2000), 123--140.

\bibitem{\refnagprob} \manyby. {\it  On Nagata's Conjecture},
preprint (2000), to appear, J. Alg.

\bibitem{\refHR} Harbourne, B. and Ro\'e, J. 
{\it Linear systems with multiple base points in \pr2}, preprint (2000).

\bibitem{\refHi} Hirschowitz, A.
{\it Une conjecture pour la cohomologie 
des diviseurs sur les surfaces rationelles g\'en\'eriques},
Journ.\ Reine Angew.\ Math. 397
(1989), 208--213.

\bibitem{\refHo} M. Homma,
{\it A souped up version of Pardini's theorem and its aplication
to funny curves}, Comp. Math. 71 (1989), 295--302.

\bibitem{\refId} Id\`a, M.
{\it The minimal free resolution for the first infinitesimal 
neighborhoods of $n$ general points in the plane},
J.\ Alg. 216 (1999), 741--753.

\bibitem{\refKa} Kac, V. G.
{\it Infinite dimensional Lie algebras},
Progress in Math. 44, Birkhauser, Boston (1983).

\bibitem{\refL} Looijenga, E.
{\it Rational surfaces with effective anticanonical divisor},
Ann. of Math. 114 (1981), 267-322.

\bibitem{\refMAC} Grayson, D., and Stillman, M.
{\it MACAULAY 2}, Version 0.8.52; archival site
{\tt www.math.uiuc.edu/Macaulay2.}

\bibitem{\refmig}
Mignon, T., {\it Syst\`emes de courbes planes \`a singularit\'es 
impos\'ees: le cas des multiplicit\'es inf\'erieures ou \'egales 
\`a quatre}, J.\ Pure Appl.\ Algebra 151 (2000), no. 2, 173--195. 

\bibitem{\refNone} Nagata, M. {\it On the 14-th problem of Hilbert}, 
Amer.\ J.\ Math.\ 33 (1959), 766--772.

\bibitem{\refNtwo} \manyby. {\it On rational surfaces, II}, 
Mem.\ Coll.\ Sci.\ 
Univ.\ Kyoto, Ser.\ A Math.\ 33 (1960), 271--293.

\bibitem{\refP} Paxia, P. {\it On flat families of fat points}, 
Proc. Amer. Math. Soc., 112 (1991).

\bibitem{\refRoe} Ro\'e, J. {\it On the existence of plane curves
with imposed multiple points}, to appear, J. Pure Appl. Alg.

\bibitem{\refRoeb} Ro\'e, J. {\it Linear systems of plane curves with
imposed multiple points}, preprint (2000).

\bibitem{\refS} Segre, B. {\it Alcune questioni su insiemi finiti
di punti in Geometria Algebrica},
Atti del Convegno Internaz. di Geom. Alg., Torino (1961).

\bibitem{\refSe} Seibert, J. {\it The Dimension of Quasi-Homogeneous 
Linear Systems With Multiplicity Four}, preprint 
({\tt http://xxx.lanl.gov/abs/math.AG/9905076}).

\bibitem{\refX} Xu, G. {\it Ample line bundles on smooth surfaces}, 
Jour.\ Reine Ang.\ Math.\ 469 (1995), 199--209.

\bye